\documentclass[12pt]{article}
\usepackage{amsmath}\usepackage{amssymb}\usepackage{amscd}
\newcounter{theorem}\makeatletter
\@addtoreset{theorem}{section}\makeatother

\newtheorem{proposition}[theorem]{Proposition}
\newtheorem{lemma}[theorem]{Lemma}\newtheorem{def-lem}[theorem]{Definition-Lemma}

\newtheorem{corollary}[theorem]{Corollary}

 \def\A{\mathcal{A}}
 \def\Ab{\bar{\A}}
 \def\Ac{\tilde{\A}}
 \def\B{\mathcal{B}}
 
 \def\be{\begin{equation}}
 \def\C{\mathbb{C}}
 \def\D{\bar{\mathcal{D}}}

 \def\der{{\bar{\partial}}}
 \def\Dif{\operatorname{Diff}}
 \def\Diffa{\operatorname{Diff}_{\h}}
 
 \def\Diff{\operatorname{\overline{Di}ff}_{\h}}
 \def\ds{\displaystyle}

 \def\ee{\end{equation}}
 
 \def\End{\operatorname{End}}
 \def\I{\operatorname{I}}
 \def\J{\operatorname{J}}
  \def\JJ{\bar{\J}}
  \def\Ker{\operatorname{Ker}}
 \def\g{{\mathbf {g}}}
 \def\gln{{\mathbf {gl}}_n}
 \def\gl{{\mathbf {gl}}}
 \def\GG{\mathcal{G}}
 \def\GGG{\bar{\mathcal{G}}}
 \def\GDif{\mathcal{G}\!\operatorname{Diff}}
 
  \def\GDiff{\mathcal{G}\!\operatorname{\overline{Di}ff}_{\h}}
 \def\h{{\mathbf{h}}}

 \def\mult{\scirc}
 \def\n{\mathbf {n}}

 \def\P{\operatorname{P}}
 \def\PP{\mathcal{P}}
 
 \def\Ph{{\mathcal{P}}_{\h}}

 \def\part{{{\partial}}}
 
 \newcommand{\rf}[1]{(\ref{#1})}
 \def\q{\check{\operatorname{q}}}

 \def\S{{\operatorname{Z}}}
 
\def\scirc{{\scriptstyle{\,\diamond\, }}}

\def\sl{{\mathbf {sl}}}
\def\t{{\mathbf t}}
\def\T{{\operatorname{T}}}

\def\th{\mathring{h}}

\def\ttDelta{\T_+}
\def\tnabla{\T_+}
\def\U{\operatorname{U}}
\def\Uh{\bar{\U}(\h)}
\def\Ug{\U(\mathbf{g})}
\def\UUg{\bar{\U}(\mathbf{g})}
\def\UU{\bar{\operatorname{U}}}
\def\ve{\varepsilon}
\def\vphi{\varphi}
\def\vf{\varphi}

\def\xii{\check{\xi}}
\def\Z{\bar{\operatorname{Z}}}
\def\Zo{\operatorname{Z}}
\def\T{\operatorname{K}}

\def\ZZ{\mathbb Z}

\begin{document}
\begin{center}
{\Large\bf Contravariant form for reduction algebras and Pieri rule}

\vspace{.4cm}
{\bf S. Khoroshkin$^\diamond$$^\star$ and
O. Ogievetsky$^{\circ}$$^\ast$\footnote{On leave of absence from P. N. Lebedev Physical Institute, Leninsky Pr. 53,
117924 Moscow, Russia} }

\vskip .2cm
$\diamond\ ${ITEP,  B.Cheremushkinskaya 25, Moscow 117218, Russia}\\
 $\star\ ${National Research University Higher School of Economics,\\  Myasnitskaya 20, Moscow 101000, Russia}\\

\vspace{.1cm}
$\circ\ ${Aix Marseille Universit\'{e}, Universit\'{e} de Toulon, \\CNRS, CPT, Marseille, France}\\
$^\ast${Kazan Federal University, Kremlevskaya 17, Kazan 420008, Russia}

\end{center}

\date{today}
\vskip .2cm
\begin{abstract}
We study properties and constructions of contravariant forms on reduction algebras. As an application we compute norms of highest weight vectors in the tensor product of an irreducible finite dimensional representation of the Lie algebra $\gln$ with a symmetric or wedge tensor power of its fundamental representation. Their zeroes describe Pieri rules.   
\end{abstract}

\section{Introduction} The contravariant (or Shapovalov) form on highest weight modules is a powerful tool in the representation theory of reductive Lie algebra. It is used for the construction of irreducible representations, description of singular vectors of Verma modules etc \cite{D}. In this paper we define an analogue of the Shapovalov form for certain reduction algebras, compute it and apply to the the well-known problems of classical representation theory, calculating the norms of $\n_+$-invariant  vectors in tensor products of irreducible finite-dimensional representation of the Lie algebra $\gln$ and symmetric or exterior powers of its fundamental representation. Here $\n_+$ is a Lie algebra of upper triangular matrices.  Zeros  of these norms describe Pieri rules.
The norms themselves can be regarded as a generalization of particular Clebsch--Gordan coefficients.

\vskip .1cm
To perform these calculations we pass to three particular reduction algebras: diagonal reduction algebra $\D(\gln)$, and $\Diff(n)$ together with its odd analogue $\GDiff(n)$, see Section \ref{secRedalgDD} for definitions. The first algebra may be regarded as a deformation of $\U(\gln)$ with coefficients in the localized universal  enveloping algebra $\Uh$, and the latter two as analogous  deformations of the algebras of polynomial differential operators in even or odd variables.
These algebras possess their own $\Uh$-valued contravariant forms whose specializations to dominant weights 
is then used for the calculation of norms of $\n_+$-invariant  vectors in tensor products of finite-dimensional representations of $\gln$.

\vskip .1cm
The main ingredient of our applications of the theory of reduction algebras 
is the connection of the contravariant form on the reduction algebra to the Zhelobenko automorphism $\xii_{w_0}$, where $w_0$ is the longest element of the Weyl group of $\gln$. The origin of this connection goes back to Zhelobenko, see \cite{Zh}. It was reformulated, proved and used in \cite{KN} 
for the description of irreducible representations of Yangians. We reproduce 
here the arguments of \cite{KN} and then compute the contravariant form on polynomial representations of the  algebras $\Diff(n)$ and $\GDiff(n)$ in two ways: using $\xii_{w_0}$ and by direct computations in the latter reduction algebras.

\vskip .1cm
The paper is organized as follows. In Sections 2.1-2.2 we recall the definition of Mickelsson algebras and their localizations called reduction algebras, introduce Zhelobenko automorphisms and describe in Section 2.3 our basic examples - reduction algebras $\D(\g)$, $\Diff(n)$ and $\GDiff(n)$. In Section 3.3 we introduce a natural class of $(\Z,\h)$-modules over reduction algebras and a notion of 
$\Uh$-valued contravariant forms on them. We establish a connection of these forms with the contravariant forms on $\n_+$-invariants and $\n_-$-coinvariants of certain $\g$-modules. Here $\g$ is a reductive Lie algebra, $\n_\pm$ are their opposite nilpotent subalgebras. In Section 3.4 we describe analogues of the Harish-Chandra map for our basic examples of reduction algebras and define with their help contravariant forms on these algebras. Section 3.5 is devoted to the calculation of these forms on basic polynomial representations of the algebras  $\Diff(n)$ and $\GDiff(n)$.   Sections 4.1-4.2 are devoted to the   justification of the evaluations of the computed contravariant forms and their use for the norms of
$\n_+$-invariant  vectors in tensor products of irreducible finite-dimensional representations of the Lie algebra $\gln$ and symmetric or exterior powers of its fundamental representation. In Section \ref{secPieri} we deduce the Pieri rules. Appendices contain an alternative derivation of norms of $\n_+$-invariant  vectors.

\section{Reduction algebras}
\subsection{Three types of reduction algebras}\label{threetypes} Let $\g$ be a finite-dimensional reductive Lie algebra with a fixed triangular decomposition  $\g=\n_++\h+\n_-$, where $\h$ is Cartan subalgebra, $\n_+$ and $\n_-$ are two opposite nilpotent subalgebras. We denote by $\Delta$ the root system of $\g$ and by $\Delta_+$ the set of positive roots.  Let $\A$ be an associative algebra which contains the universal enveloping algebra $\Ug$. In particular, $\A$ is a 
$\Ug$-bimodule with respect to the left and right multiplications by elements of $\Ug$. We assume that $\A$ is  free as the left $\Ug$-module and, moreover, that $\A$ contains a subspace $V$, invariant with respect to the adjoint action of $\Ug$ such that $\A$ is isomorphic to $\Ug\otimes V$ as the left $\Ug$ module. The action on  $\Ug\otimes V$ is diagonal. The adjoint action of $\g$ on $V$ is assumed to be reductive.

\vskip .1cm
In this setting we have three natural reduction algebras. The Mickelsson \cite{M} algebra $\S_+= \S(\A,\n_+)$ is defined as the quotient of the normalizer
of the left ideal $\J_+=\A\n_+$ modulo  $\J_+$.   The Mickelsson algebra $\S_-=\S(\A,\n_-)$ is defined as the quotient of the normalizer of the right ideal $\J_-=\n_-\A$ modulo  $\J_-$.

\vskip .1cm
In the following we assume that $\A$ is equipped with an  anti-involution $\ve$ 
whose restriction to $\Ug$ coincides with the Cartan anti-involution:
\be\label{P4}\ve(e_{\alpha_c})=e_{-\alpha_c}\ , \qquad \ve(h)=h\  \text{ for any}\ h\in\h\ ,\ee
 where   $\alpha_c$, $c=1,...,r$, are simple roots in $\Delta_+$ and $e_{\pm \alpha_c}$ and $h_{\alpha_c}=\check{\alpha}_c$ are Chevalley generators of $\g$, normalized by the conditions
 $$ [h_{\alpha_c},e_{\pm \alpha_c}]=\pm 2 e_{\pm \alpha_c},\qquad [e_{ \alpha_c}\ , e_{- \alpha_c}]= h_{\alpha_c}\ .$$
Due to \rf{P4}, 
$\ve(\J_+)=\J_-$ and $\ve(\S_+)=\S_-$
so that $\ve$ establishes an anti-isomorphism of the 
associative algebras $\S_+$ and $\S_-$

\vskip .1cm
Denote by $\T$ the multiplicative set, which consists of finite products of elements
\be\label{denom}h_\alpha+k\ ,\qquad k\in\ZZ\ .\ee 
Here $h_\alpha\in\h$ is the coroot corresponding to a root $\alpha$ of the root system $\Delta$ of the Lie algebra $\g$.
For the construction of the third reduction algebra we  localize with respect to $\T$ the enveloping algebras
$\U(\h)$, $\Ug$ and the algebra $\A$, denoting by $\UU(\h)$, $\UU(\g)$ and $\Ab$ the corresponding rings of fractions. 
Define $\Zo$ and $\Z=\Z(\A, \n_\pm)$, $\Zo\subset\Z$,  as the double coset spaces
$$\Zo=\A/(\J_-+\J_+)\ ,\ \qquad
\Z = \Ab/(\JJ_-+\JJ_+)\ ,$$
 where $\JJ_+=\Ab\n_+$ and $\JJ_-=\n_-\Ab$.
The localized double coset space $\Z$ is an associative algebra with respect to the multiplication  $\scirc$, see e.g. \cite{KO1} for details.
The multiplication $\scirc$ is described by the rule 
\be\label{mult}x\scirc y=x\P y\,\qquad \mod\JJ_++\JJ_-\ ,\ee
where $\P$ is the extremal projector \cite{AST} for $\g$, $ \P^2=\P$. The projector $\P$ belongs to a certain extension of
$\UUg$ (see \cite{Zh} for details), satisfies the properties 
\be\label{projector1}x\P=\P y=0  \ \ \ \text{for} \ \ x\in \n_+,\ y\in\n_-\ ,\ee
\be\label{projector2} \P=1 \mod \n_-\UUg,\qquad \P=1\mod
 \UUg\n_+\ ,\ee
\be\label{projector3} \ve(\P)=\P\ ,\ee
and can be given  \cite{AST} by the explicit multiplicative formula \rf{multfor}. Alternatively,
one can take representatives  $\tilde{x}\in\Ab$ and $\tilde{y}\in \Ab$ of coset classes
$x$ and $y$ such that either $\tilde{x}$ belongs to the normalizer of
the left ideal $\Ab\n_+$ or $\tilde{y}$ belongs to the normalizer of
the right ideal $\n_-\Ab$. Such representatives exist, see Lemma \ref{lemma1} (ii) below. Then $x\scirc y$ is  the image in the coset space $\Ac$
of the product $\tilde{x}\cdot\tilde{y}$. The latter description shows that the  maps
$\iota_\pm: \, \S_\pm\to\Z$, defined as compositions of natural inclusions and projections
\be\label{P1}\begin{split}
&\iota_+:\S_+=\mathrm{Norm (\J_+)}/\J_+\to \A/\J_+\to \Zo\subset \Z\ ,\\
&\iota_-:\S_-=\mathrm{Norm (\J_-)}/\J_-\to \J_-/\A\to \Zo\subset\Z\ ,
\end{split}
\ee
are homomorphisms of algebras.

\vskip .1cm
For each root $\alpha$ of the root system $\Delta$ of the Lie algebra $\g$ denote by $\th_\alpha\in\Uh$ the element
$$ \th_\alpha=h_\alpha+(\rho,h_\alpha),$$
where $\rho\in\h^*$ is the half sum of positive roots. Denote by $\ttDelta\subset \T$ the multiplicative set, which consists of finite products of elements
$(\th_\alpha+k)$ 
where $k$ is a positive integer. 
\begin{lemma}\label{lemma1}
{\slshape (i) The maps  $\iota_\pm$ are injective.

\vskip .1cm 
(ii) For each  $z\in\Zo$ there exist polynomials  $d_+,d_-\in\ttDelta$ such that $d_+\cdot z$ 
belongs to the image of $i_+$ and $z\cdot d_-$ belongs to the image of $i_-$.

\vskip .1cm  
\slshape{(iii)}  The anti-involution $\ve$ induces an anti-automorphism of the double coset algebra $\Z$, leaves invariant the subspace $\Zo$  and maps the images of ${\S}_\pm$ to 
the images of ${\S}_\mp$.
} \end{lemma}
{\it Proof.} (i) If $x\in \mathrm{Norm (\J_+)}$ then due to \rf{projector2}, $\P x\equiv x \mod \JJ_+$ (in the above mentioned extension of $\Ab$).
If $\iota_+(x)=0$ then $x\in \J_++\J_-$, but $\P \J_-=0$ by the properties of the projector thus $x\in\J_+$.

\vskip .1cm
(ii) For any $x\in\A$ the element $\P x$ (which is in the above extension of $\Ab$) belongs to the normalizer of $\J_+$ by the properties of the projectors.
Present $\P$ as a series
$\P=\sum_i d_i f_i e_i$, where $d_i$ are elements of $\Uh$, $f_i\in\U(\n_-)$, $e_i\in\U(\n_+)$. Then 
$$\P x\equiv \sum_i d_i f_i \hat{e}_i(x) \mod \JJ_+\ ,$$
where $\hat{e}_i(x)$ is the adjoint action of $e_i$ on $x$.
Since the adjoint action of $\n_+$ in $\Ab$ is locally finite, the latter sum is finite and belongs to the normalizer of $\J_+$ in $\Ab$. 
Multiplying this sum by the common multiple of $d_i$ we get the element of  $\mathrm{Norm (\J_+)}$ in $\A$. 

\vskip .1cm
(iii) Straightforward. \hfill$\square$ 

\subsection{Zhelobenko operators} 
  
It follows that the adjoint action of $\g$ on $\A$, $\hat{x}(a):=xa-ax$, $x\in\g$, $a\in\A$,
is locally finite and semisimple. That is, $\A$ can be decomposed into a direct sum of finite-dimensional $\g$-modules with respect to the adjoint action of $\g$. We assume also that simple reflections $\sigma_c$, $c=1,...,r$, of $\h$, generating the Weyl group of $\g$ are extended to automorphisms of the algebra $\A$, preserving $\Ug$. We denote them by the same symbols and assume that they still satisfy the corresponding braid group relations
\begin{equation}\label{braid}
\underbrace{\sigma_a\sigma_b\sigma_a\cdots}_{m_{ab}}=\underbrace{\sigma_b\sigma_a\sigma_a\cdots}_{m_{ab}}\ ,\qquad a,b=1,...,r, \ a\not=b\ ,
\end{equation}
where $m_{ab} = 2$ if $c_{ab} = 0$, $m_{ab} = 3$ if $c_{ab}c_{ba} = 2$ and $m_{ab}=6$ if $c_{ab}c_{ba} = 3$ 
with $c_{ab}$ the Cartan matrix of $\g$.
 
Since the adjoint action of $\g$ in $\A$ is reductive,  there is a common choice of such an extension\footnote{Other extensions by automorphisms of $\A$ 
of the  Weyl group action on $\h$ can be used here. First, one can use the inverse to \rf{braidext} or switch the positive and negative roots in \rf{braidext}.},
see e.g. \cite{K},
\begin{equation}\label{braidext}\sigma_c(x)= e^{\mathrm{ed}_{e_{\alpha_c}}}\circ e^{-\mathrm{ad}_{e_{-\alpha_c}}}\circ e^{\mathrm{ad}_{e_{\alpha_c}}}(x)\ .
\end{equation}
 
Denote by $\q_c$ the linear map $\q_c:\A\to \Ab/\JJ_+$ given by the relation

\begin{equation}\label{not7}\q_c(x):=\sum_{k\geq 0}\frac{(-1)^k}{k!}\hat{e}_{\alpha_c}^k(\sigma_c(x))
e_{-\alpha_c}^k\
\ds\prod_{j=1}^k(h_{\alpha_c}-j+1)^{-1}
\quad \mod\JJ_+\ .\end{equation} 
Properties of the operator $\q_c$ are listed in the following proposition \cite{Zh}, see also \cite{KO1}.
Here the shifted action of the Weyl group on $\h$ is used:
$$w\circ \th_\alpha=\th_{w(\alpha)}.$$
\begin{proposition} \label{prop2.2}{\slshape 
\begin{itemize} 
\item[(i)] \ $\q_c(\J_+)=0$;
\item[(ii)]\ $\q_c(\J_-)\subset \left( \JJ_-+\JJ_+\right)/\JJ_+$;
\item[(iii)]\ $\q_c(hx)=(\sigma_c\circ h)\q_c(x)$ for any $x\in\A$ and $h\in\h$;
\item[(iv)]\ $\q_c(xh)=\q_c(x)(\sigma_c\circ h)$ for any $x\in\A$ and $h\in\h$.
\end{itemize}
}\end{proposition}
The last two properties allow to extend the map $\q_c$ to the map of the localized algebras
$\q_c:\Ab\to\Ab/\JJ_+$. The properties (i) and (ii) show that the map  $\q_c$ defines a linear map of the double coset algebra $\Z$ to
itself.

\vskip .1cm
The Zhelobenko maps satisfy the braid group relations \cite{Zh}:
\begin{equation}\label{braid1}
\underbrace{\q_a\q_b\q_a\cdots}_{m_{ab}}=\underbrace{\q_b\q_a\q_a\cdots}_{m_{ab}}\ ,\qquad a,b=1,...,r, \ a\not=b
\end{equation}
and the inversion relation \cite{KO1}:
\begin{equation}\label{invr}\q_c^2(x)=(h_{\alpha_c}+1)^{-1}\ \sigma_c^2(x)\
(h_{\alpha_c}+1)\quad \mod\JJ_+\ . 
\end{equation}
In \cite{KO1} we established the following homomorphism properties of the Zhelobenko maps $\q_c$.
\begin{proposition}{\slshape
Zhelobenko map $\q_c$ defines a homomorphism of the Mickelsson algebra $\S_+$ to the double coset algebra $\Z$ and 
an automorphism of the double coset algebra $\Z$.}
\end{proposition}

One can equally start from the right ideal $\J_-$ and define Zhelobenko operators 
$ \xii_c=\ve\q_c\ve: \A\to \JJ_-\backslash\Ab$:
\begin{equation}\label{xic}\xii_c(x):=\sum_{k\geq 0}\frac{1}{k!}\prod_{j=1}^k(h_{\alpha_c}-j+1)^{-1}
{e}_{\alpha_c}^k
\hat{e}_{-\alpha_c}^k(\sigma_c(x))\
\ds
\quad \mod\JJ_-\ .\end{equation}
As well as $\q_c$ the maps $\xi_c$ determine the automorphisms 
$\xi_c:\ \Z\to\Z$ of the double coset algebra, satisfying the braid group relations \rf{braid}.  

\begin{proposition}{\slshape The following relation between automorphisms $\q_c$ and $\xii_c$ of the double coset algebra $\Z$ takes place
\begin{equation}\label{qxire}
\xii_c(x)=\q_c^{-1}\left( ( \sigma_c\ve)^2(x)\right)\ ,\end{equation}
where $x$ is a representative in $\Ab$ of the double coset. 
}\end{proposition}	
{\it Proof.} It is sufficient to check \rf{qxire} for the $\sl_2$ subalgebra $\g_c$ of $\g$ related to the simple root $\alpha_c$. The operators $\q$ and $\xii$ are automorphisms of the algebra $\Ph(2)$ (see precise definitions below) so it is sufficient to check \rf{qxire} for the 2-dimensional representation since all other representations arise as the homogeneous components of $\Ph(2)$. With the explicit formulas for $\q$, see \cite{KO3}, the calculation for the 2-dimensional representation is immediate. See also \cite{KNS}.\hfill$\square$

\vskip .1cm
Note that the automorphism  
$(\sigma_c\ve)^2$ is the involution which is $-1$ on even-dimensional irreducible representations of $\g_c$, and $+1$ on odd-dimensional irreducible representations of $\g_c$. 
 
\vskip .1cm
For $\g=\gln$, the symmetric group $S_n$ acts on the universal enveloping algebra $\Ug$ by permutation of indices. In the sequel we shall use 
this action to extend the automorphisms $\sigma_c$ of the  Weyl group action on $\h$.
 In this situation the automorphism $(\sigma_c\ve)^2$ is  identical. 
	
\subsection{Reduction algebras $\Diff(n)$ and  $\D(\g)$}\label{secRedalgDD}
In the sequel we use the following notation for the Lie algebra $\gln$. The standard generators are denoted by $e_{ij}$, the Cartan elements $e_{ii}$ by $h_i$.
We set $h_{ij}=h_i-h_j$, $\th_i=h_i-i$ and $\th_{ij}=\th_i-\th_j$. The space $\h^*$ is spanned by the elements $\epsilon_i$, $\epsilon_i(h_j)=\delta_i^j$.  

\vskip .2cm
Let $\Dif(n)$ be an associative  ring of polynomial differential operators in $n$ variables $x^{i}$, where $i=1,...,n$. 
It is generated by the elements $x^{i}$ and $\part_{i}$,  $i=1,...,n$, subject to the defining relations
\be
[x^{i},x^{j}]= [\part_{i},\part_{j}]=0\ ,\qquad [\part_{i},x^{j}]=\delta_{i}^{j}\ .
\ee
Let $\psi:\U(\gln)\to \Dif(n)$ be the homomorphism of associative algebras, such that
\be\label{psi}\psi(e_{ij})=x^{i}\part_{j}\ .\ee
Set 
\be\label{defalga}\A=\Dif(n)\otimes\U(\gln)\ .\ee
This algebra contains $\U(\gln)$ as a subalgebra generated by the elements 
$$\psi(e_{ij})\otimes 1+1\otimes e_{ij}, \qquad i,j=1,...,n\ .$$
The corresponding double coset reduction algebra $\Z$ is denoted further by $\Diff(n)$ and is 
called the algebra of $\h$-differential operators. The algebra $\Diff(n)$ is generated over $\Uh$
by the images of the elements $1\otimes x^{i}$ and $1\otimes \part_{i}$, which we denote for simplicity
by the same letters $x^{i}$ and $\part_{i}$. They satisfy quadratic relations which can be written in the $R$-matrix form, see \cite[Proposition 3.3]{KO3}.

\vskip .1cm   
As an $\Uh$-module, $\Diff(n)$ is freely generated by images in $\Z$ of elements
$1\otimes d$, where $d\in\Dif(n)$. To distinguish elements in $\Dif(n)$ and in $\Diff(n)$,
we use sometimes   the notation $:\!d\!:$ for the image in $\Diff(n)$ of a polynomial differential 
operator $d$. The anti-involution $\ve: \Dif(n)\to\Dif(n)$ is given by the rule
$$\ve(x^{i})=\part_{i}\ ,\qquad \ve(\part_{i})=x^{i}\ .$$
For the definition of the Zhelobenko operators we use the action of the symmetric group $S_n$, 
which permutes indices of the generators $x^{i}$ and $\part_{i}$. 

\vskip .1cm
The same construction applied to the ring $\GDif(n)$ of Grassmann differential operators,
generated by the odd generators $\zeta^i$ and $\delta_i$, $i=1,...,n$, with the defining relations
\begin{equation}\label{xi}
\zeta^i\zeta^j+\zeta^j\zeta^i= \delta_i\delta_j+\delta_j\delta_i=0\ ,\qquad \zeta^i\delta_j+\delta_j\zeta^i=
\delta_j^i
\end{equation}
and the homomorphism 
$\varphi:\U(\gln)\to \GDif(n)$, such that
\be\label{psi1}\varphi(e_{ij})=\zeta^{i}\delta_{j}\ee gives rise to the reduction algebra $\GDiff(n)$.
\bigskip
 
For any reductive Lie algebra $\g$  one can define the diagonal reduction algebra as follows.
Set $\A=\U(\g)\otimes\U(\g)$ and use the diagonally embedded $\U(\g)$ as $\U(\g)$-subalgebra of $\A$. 
This subalgebra is generated by the elements $x^{(1)}+x^{(2)}$, where, for $x\in\g$,
$x^{(1)}:=x\otimes 1$, 
$x^{(2)}:=1\otimes x$. The Chevalley anti-involution $\ve$ and the braid group action on $\U(\g)$ is naturally extended
to its tensor square.  The corresponding reduction algebra is denoted by $\D(\g)$ and is called the diagonal 
reduction algebra.

\vskip .1cm
There are two families of natural generators of $\D(\g)$. The first family is given by the images of
the elements $x^{(1)}$, $x\in\g$. In particular, we denote the images of Cartan-Weyl generators $e_\alpha^{(1)}$ by $s_\alpha^{(1)}$,
and the images of the elements $h_\alpha^{(1)}$, $h_\alpha\in\h$ by $t_\alpha^{(1)}$. 
    
The second family is given by the projections of the elements  $x^{(2)}$, $x\in\g$, where we use analogous notations 
with  the change of the upper index. Clearly,
$$ s_\alpha^{(1)}+ s_\alpha^{(2)}=0\ ,\qquad\text{and}\qquad t_\alpha^{(1)}+t_\alpha^{(2)}=h_\alpha\ ,\qquad \alpha\in\Delta\ .$$
   
We will be mainly interested in the diagonal reduction algebra $\D(\gln)$. The algebraic structure of the 
$\D(\gln)$ was studied in \cite{KO2,KO3}. Note that the homomorphisms \rf{psi} and \rf{psi1} define the homomorphisms of the reduction algebras
\be\label{psi2}
\psi: \D(\gln)\to \Diff(n)\qquad\text{and}\qquad \varphi:\D(\gln)\to\GDiff(n)\ .
\ee 
  
\section{ Contravariant forms} 
\subsection{Extremal projector and $\n_\mp$-(co)invariants} 
Let  $M$ be an $\A$-module. Then the space $M^\circ=M^{\n_+}$ of $\n_+$-invariants (or singular vectors or highest weight vectors)  is a $\S_+$-module, and the space 
$M_\circ=M_{\n_-}=M/\n_-M$ of $\n_-$-coinvariants is a $\S_-$-module. Assume further that $M$ is locally $\n_+$-finite, and the action of $\h$ is semisimple with {\it non-singular} (sometimes called dominant) weights, that is
$$M=\oplus_{\lambda\in\h^*} M_\lambda\ ,\qquad hv=(h,\lambda)v  \ \ v\in M_\lambda, \h\in\h\ ,$$
and 
\be\label{P3} (h_\alpha,\lambda+\rho)\not=-1,-2,...\ ,\ \alpha\in\Delta_+\ .\ee
Equivalently, the eigenvalues of all elements $\th_\alpha$, $\alpha\in\Delta_+$, are not negative integers. 
In this case, the action of the extremal projector $\P$ on $M$ is well defined, and the properties \rf{projector1} of  $\P$ 
imply that its image in $\End M$ establishes an isomorphism of $\n_+$- invariants $M^\circ$ and $\n_-$-coinvariants $M_\circ$:
\be \label{P3b}\P: M_\circ\rightsquigarrow M^\circ, \qquad x\to Px\ . \ee
If in addition the eigenvalues of all elements $\th_\alpha$, $\alpha\in\Delta_+$, are generic, that is,  
\be\label{P3a}(h_\alpha,\lambda+\rho)\not\in\ZZ\ ,\ \alpha\in\Delta_+\ ,\ee
then each of these isomorphic spaces comes equipped with a $\Z$-module structure.
The multiplication by elements  $z\in\Z$, which we sometimes denote by the symbol $\mult$ of the multiplication in the double coset algebra, can be described in several ways. First, using Lemma \ref{lemma1}, we can multiply $z$ 
by a polynomial $d_-$ from the right and get an element of $\S_-$, which we use for the action on coinvariants $M_\circ$; 
the multiplication by $d_-$ on each weight space is then given as the multiplication by a nonzero number thus is an invertible operator on $M_\circ$, so this allows to define the action of $z$ itself. 
Second, we can multiply any representative of $M_\circ$ by
$\P z\P$ (or $ z\P$  )  and get another element of $M_\circ$.
We can use analogous arguments for $M^\circ$ with the passage from $z$ to an element of $\S_+$. Finally, we can multiply an element of $M^\circ$ directly by $\P z\P$ (or $\P z$  )  and get another element of $M^\circ$. 

\vskip .2cm  
There is another special case of a natural identification of $\n_+$-invariants with $\n_-$-coinva\-ri\-ants. Assume that the restriction of an $\A$-module $M$ to $\g$ is decomposed into a direct sum of finite-dimensional $\g$-modules. In this case not all weights of $M$ are non-singular, but the weights of $M^\circ$ and of $M_\circ$ are dominant, that is, 
$$(h_\alpha,\lambda)=0,1,2,3,\dots\ ,\ \alpha\in\Delta_+\ ,$$
due to the structure of irreducible finite-dimensional $\g$-modules. Thus we have a well defined action of $\P$ on $M_\circ$ and $M^\circ$,
establishing an isomorphism of them. The  action of $\S_+$ on $M^\circ$ can be extended to the action of elements from $\Zo$, and the action of   $\S_-$ on $M_\circ$ 
can be also extended to the action of elements from $\Zo$ due to Lemma \ref{lemma1}. 

\vskip .1cm  
The functor,  attaching to a $\A$-module $M$, whose restriction  to $\g$ decomposes into a direct sum of finite-dimensional $\g$-modules, 
the $\S_+$-module $M^\circ$ is  faithful  and sends irreducible representations to irreducible representations.
To show the latter property, we choose two highest weight vectors $v,u\in M^\circ$. If $M$ is irreducible, then there exists $a\in\A$, such that
$av=u$. Then $\P a v=u$ as well.
Repeating the arguments  used in the proof of Lemma \ref{lemma1} we can replace $\P a$ by an element $d^{-1}a'$, where $a'\in\mathrm{Norm} (\J_+)$ 
and the denominator  $d\in\tnabla$ is
such that $a'v=du$. Since any highest weight of the finite dimensional module is non-singular, $d$ acts on $u$ by multiplication by a nonzero scalar $c$,
so the element $a''= c^{-1}a'\in \mathrm{Norm} (\J_+)$ maps $v$ to $u$, $a''v=u$.
An analogous picture holds for the space $M_\circ$ of coinvariants and the algebra $\S_-$.
  
\subsection{Contravariant forms} 
Let $M$ be an $\A$-module. A symmetric bilinear form $(,): M\otimes M\to \C$ is called contravariant\footnote{ The content of this section can be equally repeated for a sesquilinear contravariant form.} if
\be\label{P6a} (ax,y)=(x,\ve(a)y) \ee
for any $x,y\in M$ and $a\in\A$.
Let $M$ be an $\A$-module equipped with a contravariant form $(,)$. Then this form induces a 
pairing
$$(\ ,\,): M_\circ\otimes M^\circ\to \C$$
which is contravariant for a pair of reduction algebras $\S_-$ and $\S_+$, that is,
\be\label{P6}
(gx,y)=(x,\ve(g)y)\ ,\qquad x\in M_\circ\ ,\ y\in M^\circ\ ,\ g\in\S_-\ ,\ \ve(g)\in\S_+\ .
\ee
If $M$ is locally $\n_+$-finite, and the action of $\h$ is semisimple with non-singular weights, see \rf{P3}, then due to the isomorphism of the spaces $M_\circ$ and $M^\circ$ the contravariant form on $M$ induces the contravariant form on the space 
$M_\circ$ of $\n_-$-coinvariants, so that its value $(u,u')_\circ$ on two elements $u$ and $u'$ of $M_\circ$ is equal to
$$(u,u')_\circ= (u,\P u')\ .$$
This form satisfies the following contravariant property:
\be\label{P7}
(gu,u')=(u,\ve(g)\P u')\ee
for any $u,u'\in M_\circ$ and $g\in\S_-$.
  
On the other hand, a contravariant form on $M$ defines a symmetric bilinear form 
$(,)^\circ$ on $M^\circ$ by restriction. Under the above assumptions it satisfies the  contravariant property
\be\label{P8}
(gv,v')^\circ=(v,\P\ve(g) v')^\circ\ee
for any $u,v\in M^\circ$ and $g\in\S_+$. 
The forms on $M_\circ$ and $M^\circ$ are related as follows. For any $u,v\in M_\circ$
vectors $\P u$ and $\P u'$ belong to $M^\circ$ and
$$(u,u')=(\P u,\P u').$$
If the eigenvalues of all elements $\th_\alpha$, $\alpha\in\Delta_+$, 
in $M$ are generic, then  both forms on isomorphic spaces $M_\circ$ and $M^\circ$ are $\Z$ -contravariant,
$(gu,v)=(u,\ve(g)v)$, for any $g\in\Z$ and $u,v\in M_\circ$ (or $u,v\in M^\circ$).
  
\subsection{$(\Z,\h)$-modules}
Let now $M$ be a left module over the reduction algebra $\Z$. We call it a $(\Z,\h)$-module, or $\h$-module over the reduction algebra $\Z$ if, in addition, $M$ has a structure of a free right $\Uh$-module such that:
  
\vskip .1cm
  - $(z\mult m)\cdot h=z\mult(m\cdot h)$ for any $z\in\Z$, $m\in M$ and $h\in\Uh$;

\vskip .1cm  
  - the adjoint action of $\h$ on $M$ is semisimple.

\vskip .1cm
These conditions imply that $M$ is also free as a left $\Uh$-module. 

\vskip .1cm
For example, the reduction algebra $\Z$ itself is the $\h$-module over itself with respect to the left multiplication by elements of $\Z$ and the right multiplication by elements of $\Uh$.  

\vskip .1cm  
Assume that the weights of the adjoint action of $\h$ are generic, see \rf{P3a}. Then 
for any $\mu\in\h^*$ we can define the ``evaluation" $\Z$-module $M(\mu)$,
\be\label{P9}M(\mu)=M/MI_\mu\ee 
where $I_\mu$ is the (maximal) ideal in $\Uh$ generated by elements $h-(\mu,h)$ for all  $h\in\h$.

\vskip .1cm
We define a {\it contravariant form on an $\h$-module $M$} as a contravariant map
$(\ ,\,): M\otimes M\to \Uh$, which is  linear with respect to the right action of $\Uh$,
\be\label{formu26}\begin{split}
  	&(g\mult u,v)=(u,\ve(g)\mult v)\ ,\quad g\in\Z\ ,\\ &(uh,v)=(u,vh)=(u,v)h\ ,\quad h\in\Uh \ ,
\end{split}\ee
for any $u,v\in M$. For a generic $\mu\in\h^*$,  the evaluation 
of a contravariant form on an $\h$-module $M$ determines a $\C$-valued contravariant form on $M(\mu)$.
\medskip
  
Here is the main example, which we use in this paper, of $\h$-modules over reduction algebras.
Assume that we are given a pair $(\B,\gamma)$, which consists of an associative algebra $\B$ and an algebra homomorphism $\gamma:\U(\g)\to\B$.  Let
$\A=\B\otimes \U(\g)$. Then $\A$ contains the diagonally embedded  subalgebra $\U(\g)$, generated by the elements $\gamma(x)\otimes 1+1\otimes x$, $x\in\g$.

\vskip .1cm
Let $M$ be a $\B$-module, given as a quotient of $\B$ over its left ideal $\I$, which contains all the elements $\gamma(x)$, $x\in\n_+$. Assume that the action of 
elements of $\gamma(h)$, $h\in\h$, is semisimple and all the weights $\nu$ of this action are integers, $\nu(h_\alpha)\in\ZZ$, for any $\alpha\in\Delta$.
Consider the left $\A$ -module
$N=M\otimes (\U(\g)/\U(\g)\n_+)$. Let $\bar{N}$ be  the localization of $N$, which consists of the left fractions $d^{-1}n$, where $n\in N$ and $d$ is an element of the multiplicative set $\T$, generated by the elements $(b_\alpha+k)$, $k\in\ZZ$, where $b_\alpha=\gamma(h_\alpha)\otimes 1+1\otimes h_\alpha$ are elements of the diagonally embedded Cartan subalgebra. Define $M_{(\h)}$ to be the space of $\n_-$-coinvariants of $\bar{N}$ with respect to the diagonally embedded $\n_-$,
$$M_{(\h)}=\bar{N}/\n_-\bar{N}\ ,\ \text{where}\ N=M\otimes (\U(\g)/\U(\g)\n_+)\ .$$ 
By construction, $M_{(\h)}$ is a quotient of $\Ab$ by the sum of the right ideal $\JJ_-$ and the left ideal containing $\JJ_+$. Thus $M_{(\h)}$ is a quotient of the double coset space $\Z=\Ab/(\JJ_++\JJ_-)$ by the image in $\Z$ of some left ideal in $\Ab$. Due to the structure of the multiplication in $\Z$, $a\mult b= a\P b$, this image is also a left ideal in $\Z$, so $M_{(\h)}$ is a left $\Z$-module.  For any $m\in M$ and $h\in\h$ we set 
\begin{equation} \label{21} m\cdot h:=m (1\otimes h)\ .
\end{equation}
Since elements $1\otimes h$ normalize all the ideals defining $M_{(\h)}$, this is a well defined free right action of $\U(\h)$, commuting with the $\Z$-action on $M_{(\h)}$. Moreover, due to the integer conditions on the weights of the initial module $M$, this action has a natural extension to the action of $\Uh$,
$$ m \cdot (h_\alpha+k)^{-1}:= (h_\alpha+k-\nu(m)(h_\alpha))^{-1}m\ ,$$
where $\nu(m)\in\h^*$ is the weight of $m$. For a generic $\mu\in\h^*$ (that is, $\mu(h_\alpha)\not\in\ZZ$ for any $\alpha\in\Delta$) the specialization $M_{(\h)}(\mu)$
is isomorphic to the space of $\n_-$-coinvariants of the tensor product $M\otimes M_\mu$, where $M_\mu$ is the Verma module of $\g$ with the highest weight $\mu$,
$ M_{(\h)}(\mu)\simeq (M\otimes M_\mu)_\circ $, which is isomorphic, in its turn, to the space 
$(M\otimes M_\mu)^\circ $ of $\n_+$ - invariants, see \rf{P3b}.
Denote by $\pi_\mu: M_{(\h)}\to (M\otimes M_\mu)^\circ$ the composition of the evaluation map with the above  isomorphisms. Then $\pi_\mu(xh)=\pi_\mu(x)\mu(h)$ for any $x\in  M_{(\h)}$, $h\in\h$, and
\be\label{21a}\pi_\mu( (m\otimes 1)\cdot f)= \P \cdot (m\otimes 1_\mu)\cdot f(\mu)\ee
{for any} $ m\in M$ and $f\in\Uh$.
Here $f(\gamma)$ is the evaluation  of $f$, regarded as a 
rational function on $\h^*$ at the point $\gamma\in\h^*$. The space $(M\otimes M_\mu)^\circ $ has a natural structure of the module over the corresponding reduction algebra $\S_+$ which extends, due to conditions on $\mu$, to the structure of $\Z$ -module.

\vskip .1cm
Note also that the homomorphism $\gamma:\U(\g)\to\B$ induces the homomorphism of the reduction algebras 
\be\label{homb}
\bar{\gamma}:\D(\g)\to\Z\ee
so that $M_{(\h)}$ carries as well the structure of a $(\D(\g),\h)$-module.

\vskip .1cm
Assume that the module $M_{(\h)}$ is equipped with a contravariant form $(,)$. For a generic $\mu$ this form induces a contravariant form on $(M\otimes M_\mu)^\circ$ by the rule
\be\label{genmucompf}(\pi_\mu(x),\pi_\mu(y))=(x,y)(\mu)\ ,\ x,y\in M_{(\h)}\ .\ee
 
\vskip .2cm
The $(\Z,\h)$-modules which we use in this paper arise from the rings $\mathcal{P}(n)=\C[x^1,...,x^n]$ of polynomial functions in commuting variables and from the ring $\mathcal{G}(n)=\C[\zeta^1,...,\zeta^n]$ of polynomial functions in anti-commuting variables. The ring $\mathcal{P}(n)$ is a module over the ring $\Dif(n)$
and over the Lie algebra $\gln$. Analogously, the ring $\mathcal{G}(n)$ is a module over the ring $\mathcal{G}\!\Dif(n)$
and over the Lie algebra $\gln$. 

\vskip .1cm
We define $\Diff(n)$-module $\Ph(n)$ as a quotient of the reduction algebra
$\Diff(n)$ over the left ideal $\I_\part$, generated by all $\part_i$, $i=1,...,n$. Since 
the Cartan subalgebra normalizes the ideal $\I_\part$, $\I_\part h\subset \I_\part$ for any $h\in\h$, and the weight of any monomial is integer, we have the right action of $\Uh$ on
$\Diff(n)$ which supplies $\Ph(n)$ with a structure of $\h$-module over the reduction algebra $\Diff(n)$.
We define analogously the $\GDiff(n)$-module $\mathcal{G}_{\h}(n)$.
     
In terms of the constructions above we set $\B=\Dif(n)$, $\gamma=\psi$, see \rf{psi}, $M=\Dif(n)/\Dif(n)\{\part_1,...,\part_n\}$ in the even case and 
$\B=\GDif(n)$, $\gamma=\varphi$, see \rf{psi1}, $M=\GDif(n)/\GDif(n)\{\delta_1,...,\delta_n\}$ in the odd case.  

\vskip .2cm
{\it Example.} Let $V$ be the two-dimensional tautological representation of $\gl_2$ with the basis $\{v^1,v^2\}$. 
The $(\D\,(\gl_2),\h)$-module $V_{(\h)}$ is free as a one sided $\Uh$-module of rank 2. Its left $\D\,(\gl_2)$-module structure is described by the following formulas:
$$s_{12}^{(1)}v^1=0\ ,\ s_{12}^{(1)}v^2=v^1\frac{h_{12}}{h_{12}+1}\ ,$$
$$s_{21}^{(1)}v^1=v^2\ , s_{21}^{(1)}v^2=0\ ,$$
$$s_{11}^{(1)}v^1=v^1\ ,\ s_{11}^{(1)}v^2=v^2\frac{1}{h_{12}+1}\ ,$$
$$s_{22}^{(1)}v^1=0\ ,\ s_{22}^{(1)}v^2=v^2\frac{h_{12}}{h_{12}+1}\ ,$$
$$\th_iv^j=v^j(\th_i+\delta_i^j)\ .$$

\subsection{Harish-Chandra maps} 
Constructions of contravariant forms for reduction algebras refer to analogues of Harish-Chandra map for the universal
enveloping algebras of reductive Lie algebras. We describe this map in our two basic examples.    

\begin{lemma} \label{lemma3.1}{\slshape 
(i) \ The left ideal $\I_\part= \Diff(n)\mult\{\partial_{1},\ldots,\partial_{n}\}$ of 
$\Diff(n)$ is generated over $\Uh$  by the classes  of elements $Y\part_{i}$ where $Y\in \Dif(n)$, $i=1,...,n$. 

\vskip .1cm
(ii) \ The right ideal $\I_x= \{x^{1},\ldots,x^{n}\}\mult\Diff(n)$ of $\Diff(n)$ 
is generated over $\Uh$ by the classes  of elements $x^{i} X$ where $X\in \Dif(n)$, $i=1,...,n$. 

\vskip .1cm     
(iii) \
The natural inclusion  $\Uh\to\Diff(n)$ establishes the isomorphism of 
$\Uh$-modules $\Uh$  and $\Diff(n)/(\I_\part+\I_x)$.
}     				\end{lemma}
{\it Proof}. (i) This is a corollary of the property \rf{projector1} of the extremal projector, together with the $\mathrm{ad}_{\n_+}$-invariance of  the 
linear span of $\part_{i}$. 

\vskip .1cm
(ii) Parallel to (i). 

\vskip .1cm
(iii) Follows from the Poincar\'e--Birkhoff--Witt property of the ring $ \Diff(n)$: elements $:(x^1)^{a_1}\dots (x^n)^{a_n}
\part_1^{b_1}\dots \part_n^{b_n}:$ form a basis of $ \Diff(n)$ over $\Uh$, see \cite{KO4}. 
\hfill{$\square$}
 
\vskip .2cm    			
The map $\Diff(n)\to\Uh$, which attaches to any element $x\in\Diff(n)$ the unique 
 element $x^{(0)}\in\Uh$ such that
 $x-x^{(0)}\in \I_\part+\I_x$ is an analogue of the Harish-Chandra map $\U(\g)\to\U(\h)$. With its use we define
 in a standard way the $\Uh$-valued bilinear form on $\Diff(n)$ and on its left module 
 $\Ph(n)=\Diff(n)/\I_\part$:
 \be\label{29} (x,y)= (\ve(x)\mult y)^{(0)}.\ee
Recall that $\ve(x^i)=\partial_i$, $\ve(\partial_i)=x^i$; in particular,
$\ve(\I_\part)=\I_x$. It is not difficult to show that this form is contravariant, see \rf{formu26}, and 
symmetric 
$$ (x,y)=(y,x)$$ 
for any $x,y\in\Diff(n)$ or  $x,y\in \bar{\mathcal{P}}_{n}$. The same statements take place for $(\GDiff(n),\h)$-module $\mathcal{G}_{\h}(n)$.

\vskip .2cm
The diagonal reduction algebra $\D(\g)$ contains a family of commuting, see \cite{Zh,KO2}, elements $t_{\alpha}^{(1)}$ (in the notation of Section \ref{secRedalgDD}). Let $\C[\t]$ be the ring of polynomials in
$t_{\alpha}^{(1)}$, $\alpha\in\Delta_+$.
        
\begin{lemma} {\slshape 
(i)\ The left ideal $\I_+= \D(\g)\mult\{s_\alpha^{(1)},\alpha\in\Delta_+\}$ of 
 $\D(\g)$ is generated over $\Uh$  by the classes  of elements 
  $Y e_\alpha^{(1)}$, where $Y\in\U(\g)^{(1)}$, $\alpha\in\Delta_+$. 

\vskip .1cm     
(ii)\
The right ideal $\I_-= \{s_{-\alpha}^{(1)},\alpha\in\Delta_+\}\mult\D(\g)$ of 
$\D(\g)$ is generated over $\Uh$  by the classes  of elements 
$ e_{-\alpha}^{(1)}X$, where $X\in\U(\g)^{(1)}$, $\alpha\in\Delta_+$.

\vskip .1cm
(iii) \
The natural inclusion  $\Uh\otimes\C[\t]\to\D(\g)$ establishes the isomorphism of 
 $\Uh$-modules $\Uh\otimes\C[\t]$  and $\D(\g)/(\I_++\I_-)$.			
}\end{lemma}   
{\it Proof}. (i) This is a corollary of the property \rf{projector1} of the extremal projector, together with the $\mathrm{ad}_{\n_+}$-invariance of the linear 
space $\n_+\otimes 1$.

\vskip .1cm
(ii) Parallel to (i). 

\vskip .1cm
(iii) As in Lemma \ref{lemma3.1}, this follows from Poincar\'e--Birkhoff--Witt property of $\D(\g)$, see \cite{KO4}. \hfill{$\square$}

\vskip .2cm    			
The map $\D(\g)\to \Uh\otimes\C[\t]$, which attaches to any element $u\in\D(\g)$ the unique 
 element $u^{(0)}\in \Uh\otimes\C[\t]$ such that
 $u-u^{(0)}\in \I_++\I_-$ is an analogue of the Harish-Chandra map $\U(\g)\to\U(\h)$. With its use we define
 in a standard way the $\Uh\otimes\C[\t]$-valued bilinear form on $\D(\g)$:
 \be\label{29b} (u,v)= (\ve(u)\mult v)^{(0)}.\ee

The contravariant forms \rf{29} and \rf{29b} are compatible. Namely,  
the bilinear form \rf{29} on $\Ph(n)$ is also contravariant with respect to $\D(\gln)$ due to the homomorphism \rf{psi}.
As a $\D(\gln)$-module, $\Ph(n)$ decomposes into a direct sum of homogeneous components. 
The component $\Ph(n;k)$ of degree $k$ is generated by the element $:(x^1)^k:$ annihilated by the ideal $\I_+$. 
We have $$(:(x^1)^k:\, ,\, :(x^1)^k:)=k!\ .$$ The restriction of the form \rf{29} to $\Ph(n;k)$ can be obtained by the evaluation,
see \rf{defalga}, 
$t_{11}:=e_{11}^{(1)}\mapsto k!$, $t_{jj}:=e_{jj}^{(1)}\mapsto 0$, $j>1$, of the form \rf{29b}.
For the anti-commuting variables, the homogeneous component of degree $k$ is generated by the element 
$:\zeta^1\zeta^2\dots\zeta^k:$ annihilated by the ideal $\I_+$, for which $$(:\zeta^1\zeta^2\dots\zeta^k:\, ,\, :\zeta^1\zeta^2\dots\zeta^k:)=1\ .$$ Now the evaluation is $t_{ii}\mapsto 1$, $i=1,\dots,k$, and $t_{jj}:=e_{jj}^{(1)}\mapsto 0$, $j>k$.

\vskip .2cm
The bilinear form \rf{29} is covariant with respect to the action of Zhelobenko operators in the following sense.
\begin{lemma}\label{lemma3.3}{\slshape  For any elements $x,y\in\Ph(n)$ or $x,y\in\mathcal{G}_{\h}(n)$ we have
\be\label{covofthef}\q_c(x,y)=(\xii_c(x),\q_c(y))\ .\ee
}\end{lemma}
{\it Proof} consists in the following calculation:
$$\q_c(x,y)=\q_c((\ve(x)\mult y)^0)=(\q_c(\ve(x)\mult\q_c(y))^0=(\ve(\xii_c(x)\mult\q_c(y))^0=
(\xii_c(x),\q_c(y))\ .\qquad \square$$   

\subsection{Calculations of contravariant form on $\Ph(n)$ and $\mathcal{G}_{\h}(n)$}
We use the notation $x^{\uparrow a}=x(x+1)\dots(x+a-1)$ and $x^{\downarrow a}=x(x-1)\dots(x-a+1)$ for the Pochhammer symbols.
 
\begin{proposition} \label{prop3.3}{\slshape  Images of monomials $:x^{{\nu}}\! : = :\! (x^1)^{\nu_1}\cdots (x^n)^{\nu_n}\!: $ in  
$\Ph(n)$ have the following scalar products:
\begin{equation}\label{nuinui} 
\begin{array}{rcl}
(: x^{{\nu}} : , : x^{{\nu}'}:)&=&\delta_{{\nu},{\nu}'}\,
\displaystyle{ \prod_{k=1}^n \nu_k!\cdot \prod_{i,j:i<j}\frac{(\th_{ij}-\nu_j)^{\uparrow \nu_{i}+1}}{\th_{ij}^{\uparrow \nu_{i}+1}}}\\[2em]
&=&\displaystyle{ \delta_{{\nu},{\nu}'}\,\prod_{k=1}^n \nu_k!\cdot \prod_{i,j:i<j}\frac{\Gamma(\th_{ij}-\nu_j+\nu_i+1)\Gamma(\th_{ij})}{\Gamma(\th_{ij}-\nu_j)\Gamma(\th_{ij}+\nu_i+1)}\,.}	
\end{array}
\end{equation}    
}\end{proposition}
       
\begin{proposition} \label{prop3.4}{\slshape Images of monomials $:\zeta^{{{\nu}}}\! : = :\! (\zeta^1)^{\nu_1}\cdots (\zeta^n)^{\nu_n}\!: $ in  
$\mathcal{G}_{\h}(n)$ have the following scalar products:	
\begin{equation}
(:\zeta^{{{\nu}}}\! : , :\zeta^{{{\nu'}}}\!:)=
\delta_{{\nu},{\nu}'}\prod_{i,j:i<j}\frac{(\th_{ij}-\nu_j)^{1-\nu_i}}{(\th_{ij})^{1-\nu_i}}=
\delta_{{\nu},{\nu}'}\prod_{i,j:i<j}\frac{(\th_{ij}-1+\nu_i)^{\nu_j}}{(\th_{ij})^{\nu_j}}	
\end{equation}    
}\end{proposition}
We present two different proofs of Propositions \ref{prop3.3} and \ref{prop3.4}.

\vskip .1cm 
The first proof is based on the description of the contravariant form for certain reduction algebras given in   \cite{KN}.
We reproduce it in the particular case of the reduction algebra $\Diff(n)$.
Let $w_0$ be the longest element of the symmetric group $S_n$, regarded as the Weyl group of Lie algebra $\gln$, $w_0=(n,n-1,...,2,1)$. Let $w_0=s_{c_1}s_{c_2}\cdots  s_{c_N}$, $N=\frac{n(n-1)}{2}$, be a reduced decomposition of $w_0$. Set
$$\q_{w_0}=\q_{c_1}\q_{c_2}\cdots \q_{c_N},\qquad \xii_{w_0}=\xii_{c_1}\xii_{c_2}\cdots \xii_{c_N}.$$
Due to the braid relation \rf{braid1}, definitions of $\q_{w_0}$ and $\xii_{w_0}$ do not depend on a reduced decomposition of $w_0$ and
$$\q_{w_0} \xii_{w_0}= \xii_{w_0}\q_{w_0}=1\ .$$
For any two monomials $x^{{\nu}}=(x^1)^{\nu_1}(x^2)^{\nu_2}... (x^n)^{\nu_n}$ and 
$x^{{\nu}'}=(x^1)^{\nu'_1}(x^2)^{\nu'_2}... (x^n)^{\nu'_n}$ in commuting variables $x^1,...,x^n$, and elements $\vf_1$, $\vf_2\in\Uh$ set
\be\label{30a}<:x^{{\nu}}:\vf_1,:x^{{\nu}'}:\vf_2>=\vf_1\vf_2\;\delta_{\nu,\nu'}\;\prod_{k=1}^n \nu_k!\ .
\ee
This defines a  $\Uh$-valued pairing on $\Ph(n)$, linear with respect to the right multiplication by elements of $\Uh$. Analogously, for  any two monomials $\zeta^{{{\nu}}}=(\zeta^1)^{\nu_1}(\zeta^2)^{\nu_2}... (\zeta^n)^{\nu_n}$ and $\zeta^{{{\nu}}'}=(\zeta^1)^{\nu'_1}(\zeta^2)^{\nu'_2}... (\zeta^n)^{\nu'_n}$ in anti-commuting variables $\zeta_1,...,\zeta_n$, and elements $\vf_1$, $\vf_2\in\Uh$ set
\be\label{30b}<:\bar{\zeta}^{{{\nu}}}:\vf_1,:\bar{\zeta}^{{{\nu}}'}:\vf_2>=\vf_1\vf_2\prod_{k=1}^n \delta_{\nu_k,\nu'_k}\ .\ee
This defines a $\Uh$-valued pairing on $\mathcal{G}_{\h}(n)$, linear with respect to the right multiplication by elements of $\Uh$.
We have, see also eq. (3.18) in \cite[Section 3.3]{KN}, 
\begin{proposition}\label{prop3.5} {\slshape 
(i) For any two monomials $x^{{\nu}}$, and $x^{{\nu}'}$ the contravariant pairing of their images in $\Ph(n)$ is equal to
\be\label{30c}
(:x^{{\nu}}\!:,\, :x^{{\nu}'}\!\!:)=\q_{w_0}\left(<:w_0(x^{{\nu}}):,
\xii_{w_0}( :x^{{\nu}'}:)>\right)\ .\ee
(ii) For any two monomials $\zeta^{{{\nu}}}$, and $\zeta^{{{\nu}}'}$ the contravariant pairing of their images in $\GGG_n$ is equal to
\be\label{30d} (:\zeta^{{{\nu}}}\!:,\, :\zeta^{{{\nu}}'}\!\!:)=\q_{w_0}\left(<:w_0(\xi^{{{\nu}}}):,
 \xii_{w_0}( :\zeta^{{{\nu}}'}:)>\right)\ .\ee
 }\end{proposition}
Here in the right hand side of \rf{30c} and \rf{30d} we use the action of the symmetric group on monomials in $\PP_n$ and $\GG_n$ by permutations of indices.
The outer action of $\q_{w_0}$ is simply the shifted Weyl group action on $\Uh$. 
Proposition \ref{prop3.5} reduces the calculation of contravariant forms in $\Ph(n)$ and $\mathcal{G}_{\h}(n)$ to the calculation of Zhelobenko operator $\xii_{w_0}$, which is a simple technical exercise;
the result, e.g. for $\Ph(n)$,  is
$$ \xii_{w_0}(:x^{{\nu}}\!:)=:x^{w_0{\nu}}\!:\,\q_{w_0}\left(\prod_{i,j:i<j}
\frac{(\th_{ij}-\nu_j)^{\uparrow \nu_i+1}}{\th_{ij}^{\uparrow \nu_i+1}}\right)\ ,$$
cf \rf{nuinui}.\hfill$\square$
\medskip

\noindent{\it Remark.} The following general statement holds. Let $V$ be an irreducible finite-dimensional $\gln$-module and $<,>$ a contravariant form on $V$. Instead of \rf{30a} take its $\Uh$-linear extension to the $(\D(\gln),\h)$-module $V_{(\h)}$. 
Then the formula \rf{30c} defines a contravariant form on $V_{(\h)}$. This can be also deduced from \cite{KN}.

\vskip .2cm
\noindent{\it Proof} of Proposition \ref{prop3.5} (i). To find $(:x^{{\nu}}\!:,\, :x^{{\nu}'}\!\!:)$ we should calculate, see Lemma \ref{lemma3.1},
\be\label{kn1}\left((\ve(x^{{\nu}})\otimes 1)\Delta(\P)(x^{{\nu}'}\otimes 1)\right)^0\ee
in $\A=\Dif(n)\otimes\UU(\gln)$. Here $()^0$ means the projection of $\A$ to $\Uh$ parallel to the sum of the left ideal generated by all $\partial_i$ and diagonally embedded $e_{ij}$, $i<j$
(equivalently, by all $\partial_i$ and $1\otimes e_{ij}$, $i<j$) and of the right ideal,   generated by all $x^i$ and diagonally embedded $e_{ij}$, $i>j$
(equivalently, by all $x^i$ and $1\otimes e_{ij}$, $i>j$). The symbol $\Delta$ stands for the diagonal embedding of $\U(\gln)$. Present $\P$ in an ordered form 
 $$P=\sum_id_i(h) x_iy_i\ ,\qquad \text{where}\quad d_i\in\Uh,\ x_i\in\U(\n_-), \ y_i\in\U(\n_+)\ .$$
Then $\Delta(P)=\sum_id_i\left( h^{(1)}+h^{(2)}\right)\Delta( x_i)\Delta(y_i)$. Moving $\Delta(x_i)$ to the left and $\Delta(y_i)$ to the right, we conclude that their components in the second tensor factor do not affect the result so we can rewrite \rf{kn1} as
 \be\label{kn2}\left((\ve(x^{{\nu}}\otimes 1)) (\P[h^{(2)}]\otimes 1)(x^{{\nu}'}\otimes 1)\right)^0\ee
where $\P[h^{(2)}]$ means the shift of $\Uh$-valued coefficients in $\P$. 
The factorized expression for $\P$, see \cite{AST} for details, reads  
\be\label{multfor}
\P=\prod^{\rightarrow}_{\gamma\in\Delta_+}\P_\gamma\ \ \text{where}\ \ \P_\gamma=\sum_{n\geq 0}\frac{(-1)^n}{n! (h_\gamma+\rho(h_\gamma)+1)^{\uparrow n}}
e_{-\gamma}^ne_\gamma^n\ .\ee  
Then 
\be\notag
\P[h^{(2)}]=\prod^{\rightarrow}_{\gamma\in\Delta_+}\P_\gamma[h^{(2)}]\ \ \text{where}\ \ \P_\gamma[h^{(2)}]=
\sum_{n\geq 0}\frac{(-1)^n}{n! (h_\gamma+h^{(2)}_\gamma+\rho(h_\gamma)+1)^{\uparrow n}}e_{-\gamma}^ne_\gamma^n\ .\ee
In \rf{kn2} the elements $h_\gamma,e_{-\gamma}$ and $e_\gamma$ should be understood as differential operators, 
see \rf{psi}. Thus, the formula \rf{kn2} defines a pairing in the $\UU(\gln)$-module ${\mathcal P}(n)$ 
with coefficients in $1\otimes\Uh$, so that 
\be\label{kn4} (:x^{{\nu}}\!:,\, :x^{{\nu}'}\!\!:)=
<x^{{\nu}},\psi\left( \P[h^{(2)}]\right) x^{{\nu}'}>\ee
with the subsequent identification of elements $h^{(2)}=1\otimes h$ with elements $h\in\h$.
\smallskip
      
Next we compute  $\xii_{w_0}( :\!x^{{\nu}'}\!:)$. Since the space ${\mathcal P}(n)\otimes 1$ is an 
${\operatorname{ad}}_{\gln}$-invariant subspace of  $\A=\Dif(n)\otimes\UU(\gln)$, the consecutive application of statements (iii) and (iv) of Proposition \ref{prop2.2} leads to the following expression for $\xii_{w_0}( x^{{\nu}'})$, see \cite[Section 8.1]{KO1}:
\be\label{kn5}
\xii_{w_0}(: x^{{\nu}'}:)=\prod^{\rightarrow}_{\gamma\in\Delta_+}\sum_{n\geq 0}\frac{(-1)^n}{n!(h_{-\gamma}+\rho(h_{-\gamma})+1)^{\uparrow n}}\hat{e}^n_\gamma \hat{e}^n_{-\gamma}(w_0(x^{{\nu}'}\otimes 1))\ .\ee
Here $\hat{g}$ means as before the operator of adjoint action of $g\in\gln$. 
The adjoint action of $g\in\gln$ on ${\mathcal P}(n)$ coincides with its action 
by the left multiplication by $\psi(g)$ on ${\mathcal P}(n)$, realized as the quotient of $\Dif(n)$ over the left ideal generated by $\partial_i$. Therefore the equality \rf{kn5} can be understood
as the relation in the $\gln$-module ${\mathcal P}(n)$ with coefficients in $1\otimes \Uh$:
\be\begin{split}\notag
\xii_{w_0}(:\! x^{{\nu}'}\!:)=&\prod^{\rightarrow}_{\gamma\in\Delta_+}\sum_{n\geq 0}\frac{(-1)^n}{n!(-h^{(1)}_{\gamma}-h^{(2)}_{\gamma}-\rho(h_{\gamma})+1)^{\uparrow n}}{e}^n_\gamma {e}^n_{-\gamma}(w_0(x^{{\nu}'}))\otimes 1=\\
=&\ (w_0\otimes 1)\prod^{\rightarrow}_{\gamma\in\Delta_+}
\sum_{n\geq 0}\frac{(-1)^n}{n!(h^{(1)}_{\gamma}+h^{(2)}_{w_0(\gamma)}-\rho(h_{\gamma})+1)^{\uparrow n}}
{e}^n_{-\gamma} {e}^n_{\gamma}(x^{{\nu}'})\otimes 1\ .
\end{split}\ee	 
Here in the second line we changed the summation index from $\gamma$ to $w_0(-\gamma)$ and used the property
$\rho(h_\gamma)=-\rho(h_{w_0(\gamma)})$. We can rewrite the result using the shifted Weyl group action on the second tensor component:
\begin{align}\notag
\xii_{w_0}(:\! x^{{\nu}'}\!:)=&\ (w_0\otimes \q_{w_0})\prod^{\rightarrow}_{\gamma\in\Delta_+}\sum_{n\geq 0}\frac{(-1)^n}{n!(h^{(1)}_{\gamma}+h^{(2)}_{\gamma}+\rho(h_{\gamma})+1)^{\uparrow n}}{e}^n_{-\gamma} {e}^n_{\gamma}(x^{{\nu}'})\otimes 1  \\[.4em]
=&\ (w_0\otimes \q_{w_0})\;(\P[h^{(2)}](x^{{\nu}'})\otimes 1)\ .
\label{kn3}\end{align}
The comparison of \rf{kn2} and \rf{kn3} gives the desired statement. The proof of the part (ii) is similar.\hfill{$\square$}

\vskip .2cm
Appendix \ref{secnorms} contains the second proof of Propositions \ref{prop3.3} and \ref{prop3.4} based on explicit calculations in the rings $\Diff(n)$ and $\GDiff(n)$.

\section{Specializations}
\subsection{Specialization to non-singular weights}
The tensor product $\PP_n\otimes M_\mu$,
where $M_\mu$ is the $\gln$-Verma module with the highest weight $\mu$, is   a $\A=\Dif(n)\otimes\U(\gln)$-module, generated by the vector 
$v_\mu= 1\otimes 1_\mu$.  Since this vector satisfies the conditions
$$(\part_i\otimes 1)v_\mu= (1\otimes x) v_\mu=0\ ,\qquad i=1,...,n,\ x\in\n_+\ ,$$
and $(1\otimes h)v_\mu=\mu(h)v_\mu$ for any $h\in \h$, there is a unique $\A$-contravariant 
$\mathbb{C}$-valued form on  $\PP_n\otimes M_\mu$,
normalized by the condition $(v_\mu,v_\mu)=1$.  This contravariant form can be constructed by means of the Harish-Chandra map 
$ {}^{(0)}:\A\to \U(\h)$,
given by the prescription $x-x^{(0)}\in \I+\ve(\I)$, where $\I$ is the left ideal of $\A$, generated by $\part_i\otimes 1$, $i=1,...,n$, and $1\otimes x$, $x\in\n_+$. Then 
\be\label{31}(x\cdot v_\mu, y\cdot v_\mu)=  (\ve(x)y)^{(0)}(\mu) \qquad\text{for any}\ x,y\in\A\ .\ee
Here $(\ve(x)y)^{(0)}(\mu)$ in the right hand side of \rf{31} means the evaluation of a polynomial on elements of the Cartan subalgebra at the point $\mu\in\h^*$. The restriction of this form to the space   $(\PP_n\otimes M_\mu)^\circ$ of $\n_+$-invariants defines on this space a bilinear form, satisfying the contravariance property \rf{P8}. For a generic $\mu$, the action of the Mickelsson algebra $\Dif(n)_+=\text{Norm}(\J_+)/\J_+$ on $(\PP_n\otimes M_\mu)^\circ$  extends to the action of its localization,  the reduction algebra
$\Diff(n)$. The actions of $\Dif(n)_+$ and $\Diff(n)$ satisfy the contravariance property \rf{P6a}. 

\vskip .1cm       
Due to Lemma  \ref{lemma3.1}, the Harish-Chandra maps, defining the contravariant forms for $\Diff(n)$ and  for $\A$, are compatible, that is, they commute with the natural map from $\Dif(n)\otimes\U(\gln)$ to its double coset $\Diff(n)$; thus the contravariant form on $(\PP_n\otimes M_\mu)^\circ$ coincides with the evaluation at $\mu$ of the contravariant form on $\Ph(n)$ under the isomorphism \rf{21a}, see \rf{genmucompf}.

\vskip .1cm
We  conclude that for generic $\mu$ the square of the norm of the $\n_+$-invariant vector
$\P(x^{\nu}\otimes 1_\mu)$ of $\A$-module $\PP_n\otimes M_\mu$ is equal to, see Proposition \ref{prop3.3},
\be\label{32}\left(\P(x^{{\nu}}\otimes 1_\mu),\P(x^{{\nu}}\otimes 1_\mu)\right)=
\prod_{k=1}^n \nu_k!\cdot \prod_{i,j:i<j}\frac{(\th_{ij}(\mu)-\nu_j)^{\uparrow \nu_{i}+1}}{(\th_{ij}(\mu))^{\uparrow \nu_{i}+1}}\ .\ee
       	
On the other hand, the space $\PP_n$ decomposes into a direct sum of the spaces $S^m$ of polynomials of degree $m$, each being an irreducible $\gln$-module,
$$\PP_n=\oplus_{m\geq0}S^m.$$
The $\gln$-module $S^m$ gives rise to $(\D(\gln),\h)$-module $S^m_{(\h)}$, see Section 3.3. 
It possesses a $\D(\gln)$-contravariant form, which is the restriction of $\Diff(n)$-contravariant form on $\Ph(n)$.  
The evaluation of this form at generic $\mu$ is a $\mathbb{C}$-valued $\D(\gln)$-contravariant form on $S^m_{(\h)}$. 
Up to a normalization, the map $\pi_\mu$ transforms it to the restriction to $(S^m\otimes M_\mu)^\circ$ of the unique $\U(\gln)\otimes\U(\gln)$-contravariant form on 
$S^m\otimes M_\mu$, see \rf{genmucompf}.

\vskip .1cm 
Thus the formula \rf{32} describes norms of highest weight vectors in the tensor product of the $m$-th symmetric power 
of the fundamental representation normalized so that the square of the norm of the vector $(x^1)^m\otimes 1_\mu$ 
is equal to $m!$.    

\vskip .1cm       
Denote $\lambda=\mu+\nu$. Since the denominators of the extremal projector $\P$ belong to the set 
$\T_+$, defined in Section \ref{threetypes}, the $\n_+$-invariant vector 
$\P(x^{\nu}\otimes 1_\mu)$ is well defined for any non-singular 
$\lambda$. The square of the norm of this vector is a rational function in $\lambda$. This function 
is equal to the right hand side of \rf{32} for generic $\mu$ (that is, for generic $\lambda$). So for any non-singular $\lambda$ the right hand side 
of \rf{32} is finite and gives the square of the norm of $\P(x^{\nu}\otimes 1_\mu)$. 

\vskip .1cm
The similar considerations hold for the decomposition of the space $\mathcal{G}(n)$ into a direct sum of its homogeneous components, 
$\mathcal{G}(n)=\oplus_{m=0}^n\,\Lambda^m$ and the corresponding $(\D(\gln),\h)$-modules $\Lambda^m_{(\h)}$.
  
\subsection{Specialization to irreducible representations}
Throughout this section the weight $\lambda\in\h^*$ is assumed to be non-singular. 

\vskip .1cm
Let $M$ and $N$ be two $\g$-modules from the category $\mathcal{O}$, that is, they are $\n_+$-locally finite, $\h$-semisimple and finitely generated $\g$-modules. Assume that $M$ is generated by a highest weight vector $1_\mu$ of the weight $\mu$ (that is, M is a quotient of the Verma module $M_\mu$). Consider the $\g$-module $N\otimes M$.
\begin{lemma} \label{lemma4.1}{\slshape 
(i)\  The space $(N\otimes M)_\circ$ is spanned by the images of vectors $v\otimes 1_\mu$, $v\in N$.

\vskip .1cm		
(ii)\ For any non-singular $\lambda\in\h^*$ the space $(N\otimes M)^\circ_\lambda$ of $\n_+$-invariant 
vectors of the weight $\lambda$ is spanned by the vectors $\P(v\otimes 1_\mu)$, where $v$ has the weight $\nu=\lambda-\mu$.
}\end{lemma} 
{\it Proof}. Each element in $M$ can be presented as $g\cdot\,1_\mu$ with $g\in\U(\n_-)$. Using the relation
$$v\otimes xg\, 1_\mu\equiv -x\, v\otimes g\, 1_\mu\left(\text{mod}\ \n_-(N\otimes M)\right)\quad\text{for any } \ x\in\n_-\ , \ v\in N\ ,$$ 
we prove by induction on degree of $g\in\U(\n_-)$ that for any $v\in N$ we have an equality $v\otimes g\, v_\mu\equiv v'\otimes 1_\mu$ $\text{mod} \ \n_-(N\otimes M)$ for some $v'\in N$. This proves (i). 

\vskip .1cm
The statement (ii) follows from (i) since for a non-singular $\lambda$ the extremal projector $\P$ establishes an isomorphism of $\n_-$-coinvariants and $\n_+$-invariants of the weight $\lambda$, see Section 3.1. \hfill $\square$
   	 
Consider the tensor product $\PP_n\otimes L_\mu$ of the $\Dif(n)$-module $\PP_n$ and an irreducible 
$\U(\gln)$-module $L_\mu$ of the highest weight $\mu$ with the highest weight vector $\bar{1}_\mu$. The natural projection
$1\otimes\tau_\mu: \PP_n\otimes M_\mu\to \PP_n\otimes L_\mu$  defines  maps 
\be\label{34}\tilde{\tau}_{\mu}: (\PP_n\otimes M_\mu)^\circ\to (\PP_n\otimes L_\mu)^\circ 
\qquad\text{and}\qquad
\tilde{\tau}_{\lambda\mu}: (\PP_n\otimes M_\mu)^\circ_\lambda\to (\PP_n\otimes L_\mu)^\circ_\lambda \ee
of the spaces of $\n_+$-invariant vectors and $\n_+$-invariant vectors of the weight $\lambda$.
\begin{corollary} \label{cor4.2}(i)\ {\slshape The map $\tilde{\tau}_{\lambda\mu}$ is an epimorphism.

\vskip .1cm  	  	
(ii)\ The square of the norm of each $\n_+$-invariant vector of the weight $\lambda=\mu+\nu$ of the $\U(\gln)$-module $\PP_n\otimes L_\mu$
is given by the relation \rf{32}.
} \end{corollary}
{\it Proof}. The statement (i) follows from Lemma \ref{lemma4.1}(ii). 

\vskip .1cm
For the proof of the statement (ii) we  note that the projection map $1\otimes\tau_\mu$ 
transforms the contravariant form on $\PP_n\otimes M_\mu$ to the contravariant form on $\PP_n\otimes L_\mu$. In particular, $1\otimes\tau_\mu$  
transforms the restriction of the contravariant form to the space of $\n_+$-invariant vectors of weight $\lambda$ in $\PP_n\otimes M_\mu$ to 
the restriction of the contravariant form to the space of $\n_+$-invariant vectors of weight  $\lambda$ in $\PP_n\otimes L_\mu$,
$$(u,v)=(\tilde{\tau}_{\lambda\mu}(u), \tilde{\tau}_{\lambda\mu}(v))\ ,\ u,v\in (\PP_n\otimes L_\mu)^\circ_\lambda \ .$$
Thus for any non-singular $\lambda$ the square of the norm of the vector $\P(\bar{x}^{\nu}\otimes \bar{1}_\mu)\in 
(\PP_n\otimes L_\mu)^\circ_\lambda$
is given by the relation \rf{32}.
\hfill{$\square$} 
\bigskip
   	  
Assume now that both $\mu$ and $\lambda$ are non-singular. For any $\nu\in\h^*$ denote by $\Zo_\nu:=\Diffa(n)_\nu$ the subspace of the reduction algebra 
${\Z}=\Diff(n)$ generated by images in the double coset space $\Diff(n)$ of elements in $\A=\Dif(n)\otimes\U(\gln)$ of the weight $\nu$, 
$$\Zo_\nu:=\{ x\ \text{mod}\, (\J_++\J_-)\ \vert\ x\in\A\ ,\ [h,x]=\nu(h) x \ \ \text{for any} \ h\in\h \}\ .$$
Since the $\Dif(n)\otimes\U(\gln)$-module $\PP_n\otimes L_\mu$ is irreducible, for any vector $v\in  (\PP_n\otimes L_\mu)^\circ_\lambda$ there exists $y\in (\Dif(n)\otimes\U(\gln))_{\mu-\lambda}$ such that $y\cdot v=1\otimes \bar{1}_\mu$. Then
\be\label{35a} z\mult v:= \P y\cdot v= 1\otimes \bar{1}_\mu \ee
where $z\in \Zo_{\mu-\lambda}$ is the image of $y$ in $\Z$ . 
Due to Corollary \ref{cor4.2}, the map
$\tilde{\tau}_{\lambda\mu}: (\PP_n\otimes M_\mu)^\circ_\lambda\to (\PP_n\otimes L_\mu)^\circ_\lambda$  is an epimorphism. We now describe its kernel for dominant $\lambda$ and $\mu$ in two equivalent ways. Consider any element $u\in (\PP_n\otimes M_\mu)^\circ_\lambda$. 
\begin{lemma} \ {\slshape (i) $u\in \Ker \tilde{\tau}_{\lambda\mu}$ iff $\ z\mult u=0\ \,\text{for any}\ z\in\Zo_{\mu-\lambda } $. 	       	
\label{lemma4.2}

\vskip .1cm
(ii) $u\in \Ker \tilde{\tau}_{\lambda\mu}$ iff it is in the kernel of the contravariant form $(\,,\,)$.
}\end{lemma}
{\it Proof}. Let $u\in\Ker \tilde{\tau}_{\lambda\mu}$. Then for each $z\in\Zo_{\mu-\lambda}$ we have $z\mult u=0$. 
Indeed, the space $(\PP_n\otimes M_\mu)^\circ_\mu$ is one-dimensional and is generated by the vector $1\otimes 1_\mu$.
The map $1\otimes \tau_\mu$ is $\Dif(n)\otimes\U(\gln)$-equivariant thus the map $\tilde{\tau}_\mu$ commutes with action of elements of $\Diff(n)$. Moreover $\tilde{\tau}_\mu(1\otimes 1_\mu)=1\otimes\bar{1}_\mu$. Then the vanishing of the
left hand side of the equality
\be\label{37}\tilde{\tau}_\mu (z\mult u)= z\mult \tilde{\tau}_\mu (u)\ee
implies the relation $z\mult u=0$ since $z\mult u$ is proportional to $1\otimes 1_\mu$ by the weight reasons. On the other hand, if $z\mult u=0$ for any $z\in\Zo_{\mu-\lambda}$ then \rf{37} implies that $z\mult\tilde{\tau}_\mu u=0$ for any $z\in\Zo_{\mu-\lambda}$.
Then, by \rf{35a}, we have $ \tilde{\tau}_\mu u=0$. This proves (i). 

\vskip .1cm   	      
Next,  Lemma \ref{lemma4.1} says that each vector $v\in (\PP_n\otimes M_\mu)^\circ_\lambda$ can be presented as $x\mult (1\otimes 1_\mu)$ for some $x\in\Zo_{\lambda-\mu}$. If $v=\P (x^{\nu}\otimes 1_\mu)$ then $x=x^{\nu}\otimes 1$. Here $\nu=\lambda-\mu$. Then for each 
$u\in\Ker \tilde{\tau}_{\lambda\mu}$,
$$(u,v)=(u,x\mult(1\otimes 1_\mu))=(\ve(x)\mult u,1\otimes 1_\mu)=0$$
since $\ve(x)\in\Zo_{\mu-\lambda}$. On the other hand, if $(u,\ve(x)(1\otimes 1_\mu))=0$ for any  $x\in\Zo_{\mu-\lambda}$ then $(x\mult u, 1\otimes 1_\mu)=0$ and thus $x\mult u=0$.\hfill{$\square$} 
\begin{corollary} \label{cor4.4} {\slshape For non-singular $\lambda$ and $\mu$ the vector $\P(x^{\nu}\otimes \bar{1}_\mu)$ is 
a nonzero element of $(\PP_n\otimes L_\mu)^\circ_\lambda$ iff its norm is nonzero.
}\end{corollary}
   	      	
Let now $\mu$ be the highest weight of a finite-dimensional irreducible $\gln$-module $L_\mu$. In particular, $\mu$ is dominant. Then the weights $\lambda$ of all $\n_+$-invariant vectors of $\PP_n\otimes L_\mu$ are highest weights of finite-dimensional irreducible $\gln$-modules and are dominant; so 
they are non-singular. Corollaries \ref{cor4.2} and  \ref{cor4.4} describe all nonzero $\n_+$-invariant vectors 
of $\PP_n\otimes L_\mu$ together with their norms. 

\vskip .1cm
The considerations are valid for Grassmann variables, where now the relation	
\be\label{39}\left(\P(\zeta^{{{\nu}}}\otimes \bar{1}_\mu),\P(\zeta^{{{\nu}}}\otimes \bar{1}_\mu)\right)=
\prod_{i,j:i<j}\left( \frac{\th_{ij}(\mu)-\nu_j}{\th_{ij}(\mu)}\right)^{1-\nu_i}\ee	
describes the norms and nonvanishingness of $\n_+$-invariant vectors in the tensor product $\mathcal{G}_n\otimes L_\mu$.

\vskip .1cm
We summarize the results in the following proposition.
\begin{proposition}\label{propdom}{\slshape Assume that the weight $\mu$ is dominant.

\vskip .1cm
(i) The square of the norm of an $\n_+$-invariant vector in $\PP_n\otimes L_\mu$ is given by the relation \rf{32}.
The square of the norm of an $\n_+$-invariant vector in $\mathcal{G}_n\otimes L_\mu$ is given by the relation \rf{39}.

\vskip .1cm   	      
(ii) Any $\n_+$-invariant vector in $\PP_n\otimes L_\mu$ has a form $\P(x^{\nu}\otimes 1_\mu)$ (with dominant $\lambda=\mu+\nu$) and is nonzero iff its norm is nonzero.
}\end{proposition}

\subsection{Pieri rule}\label{secPieri} 
We recall some terminology concerning finite-dimensional representations of $\gln$. A finite-dimensional irreducible representation $L_\mu$ of highest weight 
$\mu=(\mu_1,\dots,\mu_n)$ is visualized by the Young diagram with $\mu_j$ boxes in the $j$-th row. Let $\vert\mu\vert :=\mu_1+\dots\mu_n$ and 
$\tilde{\mu}_i:=\mu_i-i$.

\vskip .1cm
The $m$-th symmetric power $L_{(m)}$  
of the tautological representation of $\gln$ corresponds to the one-row diagram with $m$ boxes. The $m$-th wedge power $L_{(1^m)}$
of the tautological representation of $\gln$ corresponds to the one-column diagram with $m$ boxes. 

\vskip .1cm
For two diagrams $\mu$ and $\lambda$, $\mu\subset\lambda$, the set-theoretical difference $\lambda\!\setminus\!\mu$ is called a horizontal strip if it contains no more than one box in any column. The difference $\lambda\!\setminus\!\mu$ is called a vertical strip if it contains no more than one box in any row.  
The Pieri rule says that for any $\mu$ the product $L_{(m)}\otimes L_\mu$ is multiplicity free and is a direct sum of $L_\lambda$ such that $\lambda\!\setminus\!\mu$ is a horizontal strip of cardinality $m$. The dual Pieri rule says that for any $\mu$ the product $L_{(1^m)}\otimes L_\mu $ is multiplicity free and is a direct sum of $V_\lambda$ such that $\lambda\!\setminus\!\mu$ is a vertical strip of cardinality $m$. The multiplicity freeness follows since the $\h$-weights of the $(\D(\gln),\h)$-modules $S^m_{(\h)}$ and $\Lambda^m_{(\h)}$ are multiplicity free.

\vskip .1cm
We keep the notation of Proposition \ref{propdom}.
Rewrite the last product in the right hand side of the formula \rf{32} in the form $\prod_{i,j:i<j}B_{i,j}^{({\nu})}(\mu)$ where 
$$B_{i,j}^{({\nu})}(\mu):=\frac{(\th_{ij}(\mu)-\nu_j)^{\uparrow \nu_i+1}}{\th_{ij}(\mu)^{\uparrow \nu_i+1}}\ .$$
Let $\lambda=\mu+\nu$ be the weight of the $\n_+$-invariant vector $\P(x^{\nu}\otimes 1_\mu)$. The denominator of $B_{i,j}^{({\nu})}(\mu)$ is positive. The 
numerator of  $B_{i,j}^{({\nu})}(\mu)$ is
$$(\tilde{\mu}_i-\tilde{\lambda}_j)(\tilde{\mu}_i-\tilde{\lambda}_j+1)\dots(\tilde{\lambda}_i-\tilde{\lambda}_j)\ .$$
The last factor is a positive integer. So the product vanishes iff 
\be\label{nerfob}\tilde{\mu}_i-\tilde{\lambda}_j\leq 0\ \text{ for some }\ i,j\ ,\ i<j\ .\ee
It is sufficient to analyze only the neighboring indices, that is, to replace \rf{nerfob} by the condition
\be\label{nerfone}\tilde{\mu}_i-\tilde{\lambda}_{i+1}\leq 0\ \text{ for some }\ i\ .\ee
Indeed, if $\tilde{\mu}_i-\tilde{\lambda}_{i+1}> 0$ and $j>i$ then $\tilde{\mu}_i-\tilde{\lambda}_{j}> 0$ since $\tilde{\lambda}_{i+1}-\tilde{\lambda}_{j}\geq 0$.
We conclude that the following lemma holds.

\begin{lemma}\label{oupie} {\slshape The tensor product $L_{(m)}\otimes L_\mu$ does not contain the representation $L_\lambda$, $\mu\subset\lambda$ and $\vert \lambda\vert =\vert\mu\vert+m$, iff the condition \rf{nerfob} holds.
}\end{lemma}
In other words, the tensor product $L_{(m)}\otimes L_\mu$ contains the representation $L_\lambda$, $\mu\subset\lambda$ and 
$\vert \lambda\vert =\vert\mu\vert+m$, if 
\be\label{nenene}\mu_i\geq\lambda_{i+1}\ \ \text{for all}\ i\ .\ee
The condition \rf{nenene} says exactly that the difference $\lambda\!\setminus\!\mu$ is a horizontal strip so we obtain the Pieri rule.

\vskip .2cm
The situation with the odd variables is different. The right hand side of \rf{39} is $\prod_{i,j:i<j}C_{i,j}^{({\nu})}(\mu)$ where 
$$C_{i,j}^{({{\nu}})}(\mu):=\left(\frac{\th_{ij}(\mu)-\nu_j}{\th_{ij}(\mu)}\right)^{1-\nu_i}\ .$$
The denominator of $C_{i,j}^{({\nu})}(\mu)$ is positive and the numerator is zero iff $\nu_i=0$ and $\nu_j=1$ for some $i,j$, $i<j$, and $\th_{ij}(\mu)=1$, or
$\mu_i-\mu_j=i-j+1$ which may occur only if $j=i+1$ and $\mu_{i+1}=\mu_i$. But then $\lambda_i=\mu_i$ and $\lambda_{i+1}=\mu_{i+1}+1=\lambda_i+1$
which cannot happen for a diagram $\lambda$. Thus all irreducible representations $L_\lambda$ such that $\nu=\lambda-\mu$ is a weight of $\Lambda^m_{(\h)}$ do appear in the tensor product
$L_{(1^m)}\otimes L_\mu$ which is exactly the statement of the dual Pieri rule about the vertical strip.

\appendix
\section{Rings $\Diff(n)$ and $\GDiff(n)$}
\subsection{$\Diff(n)$}\label{sechdiff}
The ring is generated by the elements $x^i$ and $\partial_i$.
We shall use, instead of the set of generators generators $\{x^i,\partial_i\}$ 
the set $\{x^i,\bar\der_i\}$ where, see \cite{KO4}, 
\be\label{nO8a}\bar{\der}_j= \partial_j {\vphi'_j}^{-1}\ \ \text{with}\ \ \vphi'_j=\prod_{k: k<j} \frac{\th_{jk}}{\th_{jk}-1}\ .\ee

The defining relations for 
the variables $x^i$ are 
\be\label{defvarx}x^i\scirc x^j=\frac{\th_{ij}+1}{\th_{ij}}x^j\scirc x^i\ ,\ i<j\ .\ee
The remaining defining relations read 
\begin{equation}\label{nO9a}
\begin{split}
\bar{\der}_i\scirc\bar\der_j&=\frac{\th_{ij}-1}{\th_{ij}}\,\bar\der_j\scirc \bar\der_i\ ,\qquad i<j\ ,\\
 \bar\der_j\scirc x^i&=x^i\scirc\bar\der_j, \ i>j\ ,\qquad
 \bar\der_j\scirc x^i=\frac{\th_{ij}(\th_{ij}-2)}{(\th_{ij}-1)^2}\,x^i\scirc\bar\der_j\ ,
\ i<j\ ,\\
\bar\der_i\scirc x^i &= \sum_j\frac{1}{1+\th_{ij}}\,x^j\scirc\bar\der_j+1\ .
\end{split}
\end{equation}

We have, for $i<j$: 
\be (x^i)^{\scirc a}\scirc(x^j)^{\scirc b}=(x^j)^{\scirc b}\scirc(x^i)^{\scirc a}\,\frac{(\th_{ij}+1)^{\uparrow a}}{(\th_{ij}-b+1)^{\uparrow a}}\ .\label{comrelx}\ee 
The proof is by induction, say, first on $a$ and then on $b$. 

\vskip .2cm
We have
\be x^{j_1}\scirc x^{j_2}\scirc\dots\scirc x^{j_k}=\, : x^{j_1} x^{j_2}\dots x^{j_k} \!:\ \ \text{if}\ \ j_1\geq j_2\geq\dots\geq j_k\ .\label{normalbo}\ee
The proof is by induction on $k$. Write the extremal projector in the form
$$\P=A_{2}A_{3}...A_n\ \text{ where }\ A_m=\P_{1,m}\P_{2,m}... \P_{m-1,m}\ ,\ m=2,...,n\ ,$$
with the notation $\P_{i,j}:=\P_{\epsilon_i-\epsilon_j}$, see \rf{multfor}. By the induction hypothesis, 
$x^{j_2}\scirc\dots\scirc x^{j_k}= :x^{j_2}\dots x^{j_k}:$, so
$$x^{j_1}\scirc x^{j_2}\scirc\dots\scirc x^{j_k}=x^{j_1}\scirc
: x^{j_2}\dots x^{j_k} :\, \equiv x^{j_1}\P
: x^{j_2}\dots x^{j_k} :\ .$$
The assertion \rf{normalbo} follows because $A_l : x^{j_2}\dots x^{j_k} :\, \equiv\, : x^{j_2}\dots x^{j_k} \!:$ for  $l=j_2+1,\dots,n$, and $x^{j_1}A_2...A_{j_1}\equiv x^{j_1}$.

\subsection{$\GDiff(n)$}
Now, the defining relations for 
the variables $\zeta^i$ are 
\be\label{defvarxi}\zeta^i\scirc \zeta^j=-\frac{\th_{ij}-1}{\th_{ij}}\zeta^j\scirc \zeta^i\ ,\ i<j\ .\ee
Let $\bar{\bar{\delta}}_j=\delta_j{\vphi'_j}^{-1}$.
The remaining defining relations read 
\begin{equation}\label{defrelfi}
\begin{split}
\bar{\bar{\delta}}_i\scirc\bar{\bar{\delta}}_j&=-\frac{\th_{ij}+1}{\th_{ij}}\,\bar{\bar{\delta}}_j\scirc \bar{\bar{\delta}}_i\ ,\qquad i<j\ ,\\
\bar{\bar{\delta}}_j\scirc \zeta^i&=-\zeta^i\scirc\bar{\bar{\delta}}_j\ , \ i>j\ ,\qquad
\bar{\bar{\delta}}_j\scirc \zeta^i=-\frac{\th_{ij}(\th_{ij}-2)}{(\th_{ij}-1)^2}\,\zeta^i\scirc\bar{\bar{\delta}}_j\ ,
\ i<j\ ,\\
\bar{\bar{\delta}}_i\scirc \zeta^i &=- \sum_j\frac{1}{1+\th_{ij}}\,\zeta^j\scirc\bar{\bar{\delta}}_j+1\ .
\end{split}
\end{equation}

Similarly to \rf{normalbo}, for the $\h$-Grassmann variables,
\be \zeta^{j_1}\scirc \zeta^{j_2}\scirc\dots\scirc \zeta^{j_k}=: \zeta^{j_1} \zeta^{j_2}\dots \zeta^{j_k} \!:\ \ \text{if}\ \ j_1> j_2>\dots > j_k\ .\label{normalfe}\ee

\subsection{Zhelobenko automorphisms}
Recall that we use the action of the symmetric group as the extension by the automorphisms of the Weyl group action. The action of Zhelobenko automorphisms on generators is
\be\label{nO7a}\begin{split}
\q_i(x^i)&=x^{i+1}\frac{\th_{i,i+1}}{\th_{i,i+1}-1}\ ,
\qquad \q_i(x^{i+1})=x^{i},\quad \q_i(x^j)=x^j\ ,\quad j\not=i,i+1\ ,\\
\q_i(\bar\der_{i})&=\frac{\th_{i,i+1}-1}{\th_{i,i+1}}\bar\der_{i+1}\ ,\qquad
\q_i(\bar\der_{i+1})=\bar\der_{i}\ ,\quad \q_i(\bar\der_j)=\bar\der_j\ ,\quad
j\not=i,i+1\ ,\\
\q_i(\th_j)&=\th_{s_i(j)}\ .
\end{split}
\ee
The Zhelobenko automorphisms act on $\zeta^i$ with the same coefficients as on $x^i$ .

\vskip .2cm
Let $\q$ be the Zhelobenko automorphism for $\mathbf{gl}_{2}$. Then
\be \q\bigl( (x^2)^{\scirc a}\scirc(x^1)^{\scirc b}\bigr) =(x^2)^{\scirc b}\scirc(x^1)^{\scirc a} 
\,\frac{(\th)^{\uparrow a+1}}{(\th-b)^{\uparrow a+1}}\label{actionq}\ee 
and
\be \q^{-1}\bigl( (x^2)^{\scirc a}\scirc(x^1)^{\scirc b}\bigr) =(x^2)^{\scirc b}\scirc(x^1)^{\scirc a}
\,\frac{(\th)^{\uparrow a}}{(\th-b)^{\uparrow a}}\ .\label{actionqeps}\ee 

\section{Calculation of norms}\label{secnorms}
\paragraph{Even variables.} Let ${\nu}=(\nu_1,\dots,\nu_n)$ be a milti-index, ${\nu}!=\nu_1!\dots \nu_n!$ and 
$x^{\nu} =  (x^n)^{\scirc \nu_n}\scirc\dots  \scirc(x^1)^{\scirc \nu_1}$. The monomials $x^{\nu}$ for 
${\nu}\in\mathbb{Z}_{\geq 0}^n$ form a basis of $\Ph(n)$. Define a bilinear form on $\Ph(n)$ by
\be (x^{\nu},x^{\nu'})=\delta_{{\nu},{\nu'}}{\nu}!\gamma_{{\nu}}\ \text{ where }\ \gamma_{{\nu}}:=\prod_{i,j:i<j} B_{i,j}^{({\nu})}\ \text{ and }\ 
B_{i,j}^{({\nu})}:=\frac{(\th_{ij}-\nu_j)^{\uparrow \nu_i+1}}{\th_{ij}^{\uparrow \nu_i+1}}
\ .\label{normbob}\ee
We denote by $s_i{\nu}$ the multi-index $(\nu_1,\dots,\nu_{i-1},\nu_{i+1},\nu_i,\nu_{i+2},\dots,\nu_n)$.
\begin{proposition}\label{propob1} {\slshape The form \rf{normbob} coincides with the contravariant form on $\Ph(n)$.}
\end{proposition}

By \rf{normalbo}, this proposition is equivalent to Proposition \ref{prop3.3}.

\vskip .1cm
{\it Proof} of Proposition \ref{propob1}. Since the subspaces of different $\h$-weight are orthogonal with respect to a contravariant form, 
it is sufficient to analyze the products $(x^{\nu},x^{\nu})$.
 
\vskip .1cm 
{\bf1.} We first check the covariance \rf{covofthef} of the form \rf{normbob}. Collecting pairs $B_{m,i}^{({\nu})}$ and $B_{m,i+1}^{({\nu})}$ for $m<i$ and pairs $B_{i,m}^{({\nu})}$ and 
$B_{i+1,m}^{({\nu})}$ for 
$m>i+1$ in the product for $\gamma_{{\nu}}$, we find
\be\label{xbo1}\q_i(\gamma_{{\nu}})=\gamma_{s_i{\nu}}\frac{\q_i(B_{i,i+1}^{({\nu})})}{B_{i,i+1}^{(s_i{\nu})}}=
\gamma_{s_i{\nu}}\frac{(\th_{i,i+1}+\nu_{i+1})^{\downarrow \nu_{i}+1}}{\th_{i,i+1}^{\downarrow \nu_{i}+1}}\
\frac{\th_{i,i+1}^{\uparrow \nu_{i+1}+1}}{(\th_{i,i+1}-\nu_i)^{\uparrow \nu_{i+1}+1}}\ .\ee
We have
\be\label{xbo2}\q_i(x^{\nu})=x^{s_i\nu}\;\frac{\th_{i,i+1}^{\uparrow \nu_{i+1}+1}}{(\th_{i,i+1}-\nu_i)^{\uparrow \nu_{i+1}+1}}\ \text{ and }\ 
\q_i^{-1}(x^{\nu})=x^{s_i\nu}\;\frac{(\th_{i,i+1}+\nu_{i+1})^{\downarrow \nu_{i}+1}}{\th_{i,i+1}^{\downarrow \nu_{i}+1}}\ .\ee
Therefore, the transformation laws $(x^{\nu},x^{\nu})\mapsto \left(\q_i^{-1}(x^{\nu}),\q_i(x^{\nu})\right)$  
and $(x^{\nu},x^{\nu})\mapsto \q_i(x^{\nu},x^{\nu})$ are the same so it is sufficient to prove 
\rf{normbob} for an arbitrary permutation of $(\nu_1,\dots,\nu_n)$. 

\vskip .1cm
{\bf 2.} We prove the assertion by induction on degree $\vert {\nu}\vert =\nu_1+\dots+\nu_n$, the induction base is $( 1,1)=1$. Assume that $\vert {\nu}\vert >0$. By part {\bf 1}, it is sufficient to verify the statement for ${\nu}$ such that $\nu_1>0$. 
By induction hypothesis, $(x^{\upsilon},x^{\upsilon})={\upsilon}!\gamma_{{\upsilon}}$ where ${\upsilon}=(\nu_n,\dots,\nu_2,\nu_1-1)$.
We have  
$$x^{\nu}=x^{\upsilon}\scirc x^1\ \text{ and }\ \gamma_{{\nu}}=\gamma_{{\upsilon}}\prod_{j>1}\frac{\th_{1j}-\nu_j+\nu_1}{\th_{1j}+\nu_1}\ .$$
It follows from \rf{comrelx} that 
$$x^{\upsilon}\scirc x^1=x^1\scirc x^{\upsilon}\prod_{j>1}\frac{\th_{1j}-\nu_j+\nu_1}{\th_{1j}+\nu_1}\ .$$
Therefore,
$$(x^{\upsilon}\scirc x^1,x^{\nu})=(x^1\scirc x^{\upsilon} ,x^{\nu})\prod_{j>1}\frac{\th_{1j}-\nu_j+\nu_1}{\th_{1j}+\nu_1}
=(x^{\upsilon} ,\part_1\scirc x^{\nu})\prod_{j>1}\frac{\th_{1j}-\nu_j+\nu_1}{\th_{1j}+\nu_1}\ .$$
We used the contravariance in the last equality. By (\ref{nO8a}), $\partial_1=\bar\der_1$. Now, 
$\bar\der_1\scirc x^{\nu}= (x^n)^{\scirc \nu_n}\scirc\dots \scirc (x^2)^{\scirc \nu_2}\scirc \bar\der_1\scirc (x^1)^{\scirc \nu_1}$ by (\ref{nO9a}). We have 
$$\bar\der_1\scirc x^1=1+x^1\scirc\bar\der_1+\text{linear combination of }x^j\scirc\bar\der_j\ \text{with}\ j>1\ .$$
The $\h$-derivative $\bar\der_j$, $j>1$, then moves to the right through remaining $x^1$ without a constant term, so $ \bar\der_1\scirc (x^1)^{\scirc \nu_1}=\nu_1+$ (a linear combination of terms $F_i\scirc\bar\der_i$) which does not contribute to the scalar product. \hfill$\square$

\paragraph{Odd variables.} Let now ${{\nu}}=(\nu_1,\dots,\nu_n)$ where $\nu_j\in\{0,1\}$, $j=1,\dots,n$, and  
$\zeta^{{\nu}} = (\zeta^n)^{\nu_n}\scirc\dots \scirc(\zeta^1)^{\nu_1}$. The monomials $\zeta^{{\nu}}$  
form a basis of $\mathcal{G}_{\h}(n)$. Define a bilinear form on $\mathcal{G}_{\h}(n)$ by
\be (\zeta^{{\nu}} ,\zeta^{{\nu'}} )=\delta_{{{\nu}},{{\nu'}}}\kappa_{{{\nu}}}\ \text{ where }\ \kappa_{{{\nu}}}:=\prod_{i,j:i<j} C_{i,j}^{({{\nu}})}\ \label{normfeb}\ee
and
$$C_{i,j}^{({{\nu}})}:=\left(\frac{\th_{ij}-\nu_j}{\th_{ij}}\right)^{1-\nu_i}
=\left(\frac{\th_{ij}-1}{\th_{ij}}\right)^{\nu_j(1-\nu_i)}\ .$$
We denote by $s_i{{\nu}}$ the string $(\nu_1,\dots,\nu_{i-1},\nu_{i+1},\nu_i,\nu_{i+2},\dots,\nu_n)$.
\begin{proposition} \label{propob2} {\slshape The form \rf{normfeb} coincides with the contravariant form on $\mathcal{G}_{\h}(n)$.}
\end{proposition}

By \rf{normalfe}, this proposition is equivalent to Proposition \ref{prop3.4}. The proof is along the same lines as for Proposition \ref{propob1}. 

\vskip .1cm
{\it Proof.} {\bf1.} Analogues of formulas \rf{xbo1} and \rf{xbo2} are
$$\q_i(\kappa_{{{\nu}}})=\kappa_{s_i{{{\nu}}}}
\left(\frac{\th_{i,i+1}+\nu_{i+1}}{\th_{i,i+1}}\right)^{1-\nu_i}
\left(\frac{\th_{i,i+1}}{\th_{i,i+1}-\nu_i}\right)^{1-\nu_{i+1}}\ ,$$
$$\q_i(\zeta^{{\nu}})\!=\!(-1)^{\nu_i\nu_{i+1}}\zeta^{s_i\boldsymbol \nu}\left(\!\frac{\th_{i,i+1}}{\th_{i,i+1}-\nu_i}\!\right)^{1-\nu_{i+1}}\!\!,\ \q_i^{-1}(\zeta^{{\nu}})\!=\!(-1)^{\nu_i\nu_{i+1}}\zeta^{s_i\boldsymbol \nu}
\left(\!\frac{\th_{i,i+1}+\nu_{i+1}}{\th_{i,i+1}}\!\right)^{1-\nu_i}\ .$$
Again, the transformation laws $(\zeta^{{\nu}},\zeta^{{\nu}})\mapsto 
\left(\q_i^{-1}(\zeta^{{\nu}}),\q_i(\zeta^{{\nu}})\right)$  
and $(\zeta^{{\nu}},\zeta^{{\nu}})\mapsto \q_i(\zeta^{{\nu}},\zeta^{{\nu}})$ 
are the same so it is sufficient to prove \rf{normfeb} for an arbitrary permutation of $(\nu_1,\dots,\nu_n)$. 

\vskip .1cm
{\bf 2.} Induction now is on degree $\vert {{\nu}}\vert =\nu_1+\dots+\nu_n$. Assume that 
$\vert {{\nu}}\vert >0$. By part {\bf 1}, it is sufficient to verify the statement for ${{\nu}}$ 
such that $\nu_1=1$. 
By induction hypothesis, $(\zeta^{\upsilon},\zeta^{\upsilon})=\kappa_{{\upsilon}}$ where 
${\upsilon}=(\nu_n,\dots,\nu_2,0)$.
We have  
$$\zeta^{{\nu}}=\zeta^{\upsilon}\scirc \zeta^1\ \text{ and }\ \kappa_{{{\nu}}}=
\kappa_{{\upsilon}}\prod_{j>1}\frac{\th_{1j}}{\th_{1j}-÷\nu_j}\ .$$
On the other hand,
$$\zeta^{\upsilon}\scirc \zeta^1=(-1)^{\vert {\upsilon}\vert }\zeta^1\scirc \zeta^{\upsilon}\prod_{j>1}\left(\frac{\th_{1j}}{\th_{1j}-1}\right)^{\nu_j}=
(-1)^{\vert {\upsilon}\vert }\zeta^1\scirc \zeta^{\upsilon}\prod_{j>1}\frac{\th_{1j}}{\th_{1j}-\nu_j}\ .$$
The rest of the proof follows, as for the even variables, from the covariance and the fact that $\bar{\bar{\delta}}_1=\delta_1$.
\hfill$\square$
 
\vskip .2cm
\noindent{\bf Acknowledgments.} We thank CIRM, Marseille, where a part of this work was done, for the hospitality.  The work of both authors was supported by the grant RFBR 17-01-00585. The work of O. O. was also supported by the Program of Competitive Growth of Kazan Federal University and the work of S. K. by the  Russian Academic Excellence Project ``5-100".

\end{document}